\documentclass{amsart}

\usepackage{amsmath}
\usepackage{amsfonts}
\usepackage{amssymb}

\numberwithin{equation}{section}

\theoremstyle{plain}

\newtheorem{theorem}[subsection]{Theorem}
\newtheorem{proposition}[subsection]{Proposition}
\newtheorem{corollary}[subsection]{Corollary}
\newtheorem{lemma}[subsection]{Lemma}

\theoremstyle{definition}

\newtheorem{conjecture}[subsection]{Conjecture}

\newtheorem{question}[subsection]{Question}

\providecommand{\supp}{\mathop{\rm supp}\nolimits}
\providecommand{\codim}{\mathop{\rm codim}\nolimits}

\newcommand{\wh}{\widehat}
\newcommand{\wt}{\widetilde}

\begin{document}

\title{The Littlewood-Gowers problem}

\author{T. Sanders}

\address{Department of Pure Mathematics and Mathematical Statistics\\
Centre for Mathematical Sciences\\
Wilberforce Road\\
Cambridge CB3 0WA\\
England } \email{t.sanders@dpmms.cam.ac.uk}

\begin{abstract} The paper has two main parts. To begin with suppose that $G$ is a compact Abelian group.
Chang's Theorem can be viewed as a structural refinement of Bessel's
inequality for functions $f \in L^2(G)$. We prove an analogous
result for functions $f \in A(G)$, where $A(G)$ is the space $\{f
\in L^1(G): \|\wh{f}\|_1 < \infty\}$ endowed with the norm
$\|f\|_{A(G)}:=\|\wh{f}\|_1$, and generalize this to the approximate
Fourier transform on Bohr sets.

As an application of the first part of the paper we improve a recent
result of Green and Konyagin: Suppose that $p$ is a prime number and
$A \subset \mathbb{Z}/p\mathbb{Z}$ has density bounded away from 0
and 1 by an absolute constant. Green and Konyagin have shown that
$\|\chi_A\|_{A(\mathbb{Z}/p\mathbb{Z})} \gg_{\varepsilon} (\log
p)^{1/3-\varepsilon}$; we improve this to
$\|\chi_A\|_{A(\mathbb{Z}/p\mathbb{Z})} \gg_{\varepsilon} (\log
p)^{1/2-\varepsilon}$. To put this in context it is easy to see that
if $A$ is an arithmetic progression then
$\|\chi_A\|_{A(\mathbb{Z}/p\mathbb{Z})} \ll \log p$.
\end{abstract}

\maketitle

\section{Introduction}

We use the Fourier transform on compact Abelian groups, the basics
of which may be found in Chapter 1 of Rudin \cite{WR}; we take a
moment to standardize our notation.

Suppose that $G$ is a compact Abelian group. Write $\wh{G}$ for the
dual group, that is the discrete Abelian group of continuous
homomorphisms $\gamma:G \rightarrow S^1$, where $S^1:=\{z \in
\mathbb{C}:|z|=1\}$. Although the natural group operation on
$\wh{G}$ corresponds to pointwise multiplication of characters we
shall denote it by `$+$' in alignment with contemporary work. $G$
may be endowed with Haar measure $\mu_G$ normalised so that
$\mu_G(G)=1$ and as a consequence we may define the Fourier
transform $\wh{.}:L^1(G) \rightarrow \ell^\infty(\wh{G})$ which
takes $f \in L^1(G)$ to
\begin{equation*}
\wh{f}: \wh{G} \rightarrow \mathbb{C}; \gamma \mapsto \int_{x \in
G}{f(x)\overline{\gamma(x)}d\mu_G(x)}.
\end{equation*}
We write
\begin{equation*}
A(G):=\{f \in L^1(G): \|\wh{f}\|_1 < \infty\}
\end{equation*}
and define a norm on $A(G)$ by $\|f\|_{A(G)}:=\|\wh{f}\|_1$.

We are concerned with the following analogue of a question due to
Littlewood, closely related to a problem of Gowers', which was
first addressed by Green and Konyagin in \cite{BJGSVK}.
\begin{question}
Suppose that $p$ is a prime number and $A \subset
\mathbb{Z}/p\mathbb{Z}$ has density bounded away from 0 and 1 by
an absolute constant. How small can
$\|\chi_A\|_{A(\mathbb{Z}/p\mathbb{Z})}$ be?
\end{question}
In their paper Green and Konyagin prove the following result.
\begin{theorem}\label{BJGSVKlpthm}
Suppose that $p$ is a prime number and $A \subset
\mathbb{Z}/p\mathbb{Z}$ has density bounded away from 0 and 1 by
an absolute constant. Then
\begin{equation*}
\|\chi_A\|_{A(\mathbb{Z}/p\mathbb{Z})} \gg \left(\frac{\log
p}{\log \log p}\right)^{1/3}.
\end{equation*}
\end{theorem}
By analogy with the original problem of Littlewood they observe
that more is probably true, indeed one might make the following
conjecture.
\begin{conjecture}\label{bkc}\emph{(Green-Konyagin-Littlewood conjecture)}
Suppose that $p$ is a prime number and $A \subset
\mathbb{Z}/p\mathbb{Z}$ has density bounded away from 0 and 1 by
an absolute constant. Then
\begin{equation*}
\|\chi_A\|_{A(\mathbb{Z}/p\mathbb{Z})} \gg \log p.
\end{equation*}
\end{conjecture}
Certainly no more than this is true as any arithmetic progression
of density bounded away from 0 and 1 shows.

In this paper we improve Theorem \ref{BJGSVKlpthm}, increasing the
exponent of $\log p$ from $1/3-\varepsilon$ to $1/2-\varepsilon$.
Specifically we show the following.
\begin{theorem}\label{maintheorem}
Suppose that $p$ is a prime number and $A \subset
\mathbb{Z}/p\mathbb{Z}$ has density bounded away from 0 and 1 by
an absolute constant. Then
\begin{equation*}
\|\chi_A\|_{A(\mathbb{Z}/p\mathbb{Z})} \gg \left(\frac{\log p}{(\log
\log p)^3}\right)^{1/2}.
\end{equation*}
\end{theorem}
The rest of the paper breaks into three parts as follows.
\begin{itemize}
\item \S\S\ref{struct}-\ref{remstruct}: In these sections we
develop a new structure theorem (Theorem \ref{AGChang}) for the
Fourier spectrum of functions in $A(G)$. \item
\S\S\ref{localft}-\ref{localstruct}: Here we recall the basic
facts of local Fourier analysis and localize our new structure
theorem in Proposition \ref{localagchang}. \item
\S\S\ref{tfr}-\ref{comp}: Finally we address the problem in this
introduction and use the local version of our structure theorem to
prove Theorem \ref{maintheorem}.
\end{itemize}
It is the results of the first two parts which seem most likely to
have further uses; Theorem \ref{maintheorem} should be seen as an
application, albeit the motivating one, of these new tools.

\section{A structural result for the Fourier spectrum}\label{struct}

One might begin studying the structure of the Fourier spectrum of
functions by looking at the sets of characters at which $\wh{f}$
is large; a natural realization of these sets is the sets
\begin{equation*}
\{\gamma \in \wh{G}: |\wh{f}(\gamma)| \geq \epsilon \|f\|_1\}
\textrm{ for } \epsilon \in (0,1].
\end{equation*}
The study of these sets has been surveyed by Green in \cite{BJGMS}
so we are brief and only recall the key facts. Write $\Gamma= \{
\gamma \in \wh{G}: |\wh{f}(\gamma)| \geq \epsilon \|f\|_1\}$.
Bessel's inequality yields
\begin{equation}\label{bessel}
|\Gamma|\epsilon^2\|f\|_1^2 \leq  \|\wh{f}\|_2^2 \leq \|f\|_2^2
\Rightarrow |\Gamma| \leq \epsilon^{-2}(\|f\|_2\|f\|_1^{-1})^2.
\end{equation}
Since $G$ is compact the norms nest and in particular
$\|f\|_2\|f\|_1^{-1} \geq 1$. Now, there is a result of Chang from
\cite{MCC} which refines (\ref{bessel}) if $\|f\|_2\|f\|_1^{-1}$ is
much larger than 1. We require some further notation to state this.
If $\Lambda$ is a set of characters on $G$ and $m \in
\mathbb{Z}^\Lambda$ has finite support then put
\begin{equation*}
m.\Lambda:=\sum_{\lambda \in \Lambda}{m_\lambda.\lambda} \textrm{
and } |m|:=\sum_{\lambda \in \Lambda}{|m_\lambda|},
\end{equation*}
where the second `$.$' in the first definition is the natural action
of $\mathbb{Z}$ on $\wh{G}$. Write $\langle \Lambda \rangle$ for the
set of all finite $\pm$-sums of elements of $\Lambda$, that is
\begin{equation*}
\langle \Lambda \rangle := \left\{m.\Lambda: m \in
\{-1,0,1\}^\Lambda \textrm{ and }|m| < \infty\right\}.
\end{equation*}
\begin{theorem}\label{chang}
\emph{(Chang's Theorem)} Suppose that $G$ is a compact Abelian
group, $f \in L^2(G)$ and $\Gamma=\{\gamma \in \wh{G}:
|\wh{f}(\gamma)| \geq \epsilon \|f\|_1\}$ for some $\epsilon \in
(0,1]$. Then there is a set of characters $\Lambda$ such that
$\Gamma \subset \langle \Lambda \rangle$ and
\begin{equation*}
|\Lambda| \ll \epsilon^{-2} (1+\log \|f\|_2\|f\|_1^{-1} ).
\end{equation*}
\end{theorem}
We develop an analogue of Chang's Theorem with $A(G)$ in place of
$L^2(G)$; here it turns out that the natural realization of the
sets of characters at which $\wh{f}$ is large is the sets
\begin{equation*}
\{\gamma \in \wh{G}: |\wh{f}(\gamma)| \geq \epsilon \|f\|_\infty\}
\textrm{ for }\epsilon \in (0,1].
\end{equation*}
There is an easy analogue of (\ref{bessel}): Write $\Gamma= \{
\gamma \in \wh{G}: |\wh{f}(\gamma)| \geq \epsilon \|f\|_\infty\}$.
Then
\begin{equation}\label{l1bessel}
|\Gamma|\epsilon\|f\|_\infty \leq  \|\wh{f}\|_1 = \|f\|_{A(G)}
\Rightarrow |\Gamma| \leq
\epsilon^{-1}\left(\|f\|_{A(G)}\|f\|_\infty^{-1}\right).
\end{equation}
A trivial instance of Hausdorff's inequality tells us that
$\|f\|_{A(G)}\|f\|_\infty^{-1} \geq 1$, and indeed the quantity
$\|f\|_{A(G)}\|f\|_\infty^{-1}$ plays the same r\^{o}le in $A(G)$
as the quantity $\|f\|_2\|f\|_1^{-1}$ does in $L^2(G)$. To
complete the square then, we shall prove the following.
\begin{theorem}\label{AGChang}
Suppose that $G$ is a compact Abelian group, $f \in A(G)$ and
$\Gamma=\{\gamma \in \wh{G}:|\wh{f}(\gamma)| \geq \epsilon
\|f\|_\infty\}$ for some $\epsilon \in (0,1]$. Then there is a set
of characters $\Lambda$ such that $\Gamma \subset \langle \Lambda
\rangle$ and
\begin{equation*}
|\Lambda| \ll \epsilon^{-1} (1+\log \|f\|_{A(G)}\|f\|_\infty^{-1}).
\end{equation*}
\end{theorem}

\section{The proof of Theorem \ref{AGChang}}\label{proofagchang}

We say that a set of characters $\Lambda$ is \emph{dissociated} if
\begin{equation*}
m \in \{-1,0,1\}^{\Lambda} \textrm{ and } m.\Lambda=0_{\wh{G}}
\textrm{ imply that } m \equiv 0.
\end{equation*}
We have the following simple lemma regarding dissociated sets.
\begin{lemma}\label{maxdis}
Suppose that $G$ is a compact Abelian group, $\Gamma$ is a set of
characters on $G$ and $\Lambda$ is a maximal dissociated subset of
$\Gamma$. Then $\Gamma \subset \langle \Lambda \rangle$.
\end{lemma}
To prove this one supposes, for a contradiction, that there is a
$\gamma \in \Gamma\setminus\langle \Lambda \rangle$. If one adds
this $\gamma$ to $\Lambda$ it is easy to see that the resulting
set is strictly larger and dissociated.

In view of this lemma Theorem \ref{AGChang} follows from:
\begin{proposition}\label{dissprop}
Suppose that $G$ is a compact Abelian group, $f \in A(G)$, and
$\Lambda$ a dissociated subset of $\{\gamma \in \wh{G}:
|\wh{f}(\gamma)| \geq \epsilon \|f\|_\infty\}$ for some $\epsilon
\in (0,1]$. Then
\begin{equation*}
|\Lambda| \ll \epsilon^{-1} (1+\log \|f\|_{A(G)}\|f\|_\infty^{-1}).
\end{equation*}
\end{proposition}
We prove this using a standard inner product type argument for
which we require an auxiliary measure.
\begin{proposition}\label{auxiliary} \emph{(Auxiliary measure)}
Suppose that $G$ is a compact Abelian group, $\Lambda$ a finite
dissociated set of characters on $G$ and $\omega \in
\ell^\infty(\Lambda)$ has $\|\omega\|_\infty \leq 1$. Then for any
$\eta \in (0,1]$ there is a measure $\mu_\eta \in M(G)$ such that
\begin{equation*}
\wh{\mu_\eta}|_{\Lambda}=\omega, \|\mu_{\eta}\| \ll (1+\log
\eta^{-1}) \textrm{ and } |\wh{\mu_{\eta}}(\gamma)| \leq \eta
\textrm{ for all } \gamma \not \in \Lambda.
\end{equation*}
\end{proposition}
\begin{proof}
[Proof of Proposition \ref{dissprop}] We define
\begin{equation*}
\omega(\lambda):=\frac{\wh{f}(\lambda)}{|\wh{f}(\lambda)|}
\textrm{ for all } \lambda \in \Lambda.
\end{equation*}
$\omega \in \ell^\infty(\Lambda)$ and $\|\omega\|_\infty \leq 1$
so we may apply Proposition \ref{auxiliary} to get the auxiliary
measure $\mu_\eta$ corresponding to $\omega$. We examine the inner
product $\langle f,\mu_\eta \rangle$.
\begin{eqnarray*}
|\langle f,\mu_\eta \rangle| & = & |\sum_{\gamma \in
\wh{G}}{\wh{f}(\gamma)\overline{\wh{\mu_\eta}(\gamma)}}|\textrm{ by Plancherel's theorem,}\\
& = & |\sum_{\lambda \in
\Lambda}{\wh{f}(\lambda)\overline{\wh{\mu_\eta}(\lambda)}}+
\sum_{\gamma \not \in
\Lambda}{\wh{f}(\gamma)\overline{\wh{\mu_\eta}(\gamma)}}|\\ & \geq &
|\sum_{\lambda \in
\Lambda}{\wh{f}(\lambda)\overline{\wh{\mu_\eta}(\lambda)}}| -
|\sum_{\gamma \not \in
\Lambda}{\wh{f}(\gamma)\overline{\wh{\mu_\eta}(\gamma)}}|\\ &\geq &
|\sum_{\lambda \in
\Lambda}{\wh{f}(\lambda)\overline{\omega(\lambda)}}| -
\eta\sum_{\gamma \not \in \Lambda}{|\wh{f}(\gamma)|}\textrm{ from the properties of $\mu_\eta$,}\\
& \geq &
\sum_{\lambda \in \Lambda}{|\wh{f}(\lambda)|} - \eta\|f\|_{A(G)}\\
& \geq & |\Lambda|\epsilon \|f\|_\infty - \eta\|f\|_{A(G)}.
\end{eqnarray*}
However
\begin{equation*}
|\langle f,\mu_\eta \rangle| \leq \|f\|_\infty \|\mu_\eta\|
 \ll \|f\|_\infty (1+\log \eta^{-1}),
\end{equation*}
so that
\begin{equation*}
\|f\|_\infty (1+\log  \eta^{-1})\gg |\Lambda|\epsilon \|f\|_\infty -
\eta\|f\|_{A(G)}.
\end{equation*}
Choosing $\eta^{-1}=\|f\|_{A(G)}\|f\|_\infty^{-1}$ yields the
result.
\end{proof}

\section{Constructing the auxiliary
measure}\label{auxiliarymeasure}

The construction of the auxiliary measure is best illustrated in
the model setting of $\mathbb{F}_2^n$ where we benefit from two
simplifications. Suppose that $\Lambda$ is a set of characters on
$\mathbb{F}_2^n$. Then
\begin{itemize}
\item $\langle \Lambda \rangle$ is simply the subspace of $\wh{G}$
generated by $\Lambda$; \item $\Lambda$ is dissociated if and only
if it is linearly independent over $\wh{\mathbb{F}_2^n}$.
\end{itemize}
The first of these is simply a convenience while the second
represents the major obstacle in transferring the arguments of
this section to the general setting. We shall prove the following
result.
\begin{proposition}\label{Mela}
Suppose that $\Lambda$ is a linearly independent set of characters
on $\mathbb{F}_2^n$ and $\omega :\Lambda \rightarrow [-1,1]$. Then
for any $\eta \in (0,1]$ there is a measure $\mu_\eta \in
M(\mathbb{F}_2^n)$ such that
\begin{equation*}
\wh{\mu_\eta}|_{\Lambda}=\omega, \|\mu_{\eta}\| \ll (1+\log
\eta^{-1}) \textrm{ and } |\wh{\mu_{\eta}}(\gamma)| \leq \eta
\textrm{ for all } \gamma \not \in \Lambda.
\end{equation*}
\end{proposition}
We relegate the technical process of extending this construction
to arbitrary compact Abelian groups to Appendix \ref{app1}.

Riesz products are the building blocks of the measure, the
necessary details of which we now record.

\subsection{Riesz products}\label{RPs}

Suppose that $\Lambda$ is a finite set of characters. If
$\omega:\Lambda \rightarrow [-1,1]$ then we define the product
\begin{equation}\label{rieszproduct}
p_\omega:=\prod_{\lambda \in
\Lambda}{\left(1+\omega(\lambda)\lambda\right)}.
\end{equation}
Such a product is called a \emph{Riesz product} and it is easy to
see that it is real and non-negative from which it follows that
$\|p_{\omega}\|_1 = \wh{p_{\omega}}(0_{\wh{G}})$. Further,
expanding out the product reveals that, $\supp \wh{p_\omega}
\subset \langle \Lambda \rangle$.

If $\Lambda$ is linearly independent then we can easily compute
the Fourier transform of a Riesz product. Suppose that $\gamma \in
\langle \Lambda \rangle$, then there is a unique $m:\Lambda
\rightarrow \{0,1\}$ such that $\gamma=m.\Lambda$ by the linear
independence of $\Lambda$, so
\begin{equation*}
\wh{p_{\omega}}(\gamma):=\prod_{{\lambda \in \Lambda} \atop
{m_\lambda \neq 0}}{\omega(\lambda)}.
\end{equation*}
This leads to the observation that
$\|p_\omega\|_1=\wh{p_{\omega}}(0_{\wh{G}})=1$ and
$\wh{p_\omega}|_\Lambda = \omega$. Moreover if $t \in [-1,1]$ then
\begin{equation*}
\wh{p_{t\omega}}(m.\Lambda):=t^{|m|}\wh{p_{\omega}}(m.\Lambda)
\textrm{ where, as before, } |m|=\sum_{\lambda \in
\Lambda}{|m_\lambda|},
\end{equation*}
So, if $|m|>1$ then
\begin{equation*}
|\wh{p_{t\omega}}(m.\Lambda)| \leq
|t|^2|\wh{p_{\omega}}(m.\Lambda)| \leq |t|^2\|p_\omega\|_1 = |t|^2
\end{equation*}
It follows that a lot of the Fourier coefficients of $p_{t\omega}$
are already small if $|t|$ is small. By taking
$\mu_\eta:=\eta^{-1}p_{\eta \omega}$ we get a well known primitive
version of the auxiliary measure of Proposition \ref{auxiliary}.
\begin{proposition}\label{primitiveauxiliary} \emph{(Primitive auxiliary measure)}
Suppose that $\Lambda$ is a linearly independent set of characters
on $\mathbb{F}_2^n$ and $\omega:\Lambda \rightarrow [-1,1]$. Then
for any $\eta \in (0,1]$ there is a measure $\mu_\eta \in
M(\mathbb{F}_2^n)$ such that
\begin{equation*}
\wh{\mu_\eta}|_{\Lambda}=\omega, \|\mu_{\eta}\| \ll \eta^{-1}
\textrm{ and } |\wh{\mu_{\eta}}(\gamma)| \leq \eta \textrm{ for
all } \gamma \not \in \Lambda \cup \{0_{\wh{G}}\}.
\end{equation*}
\end{proposition}
The basic idea for improving the measure of Proposition
\ref{primitiveauxiliary} rests on the observation that if $|m|$ is
large then $|\wh{p_{t\omega}}(m.\Lambda)|$ is in fact guaranteed to
be \emph{very} small. To construct a better measure we take linear
combinations of Riesz products so that their Fourier transforms
cancel on the characters $m.\Lambda$ where $|m|$ is small (except of
course for $|m|=1$). Begin by considering
\begin{equation*}
\nu_t := \frac{1}{2}(p_{t\omega}-p_{-t\omega}).
\end{equation*}
Then
\begin{equation*}
\wh{\nu_t}|_\Lambda = t\omega, \|\nu_t\| \leq
1,|\wh{\nu_t}(m.\Lambda)| \leq t^{|m|}
\end{equation*}
and
\begin{equation*}
\wh{\nu_t}(m.\Lambda)=0 \textrm{ if } |m| \equiv 0 \pmod 2.
\end{equation*}
It follows that
\begin{equation*}
\wh{\nu_t}|_{\Lambda}=t\omega, \|\nu_t\| \leq 1 \textrm{ and }
|\wh{\nu_t}(\gamma)| \leq t^3 \textrm{ for all } \gamma \not \in
\Lambda.
\end{equation*}
If we put $\mu_\eta=\nu_{\sqrt{\eta}}$ then we have a refinement
of Proposition \ref{primitiveauxiliary} with $\|\mu_\eta\| \ll
\eta^{-1/2}$ instead of $\|\mu_\eta\| \ll \eta^{-1}$.

More generally we consider a measure $\tau$ on $[-1,1]$ and put
\begin{equation*}
\nu_\tau:=\int{p_{t\omega}d\tau(t)}.
\end{equation*}
Then
\begin{equation*}
\|\nu_\tau\| \leq \sup_{t \in [-1,1]}{\|p_{t\omega}\|}.\|\tau\| =
\|\tau\| \textrm{ and } \wh{\nu_\tau}(m.\Lambda) =
\int{t^{|m|}d\tau(t)}\wh{p_{\omega}}(m.\Lambda).
\end{equation*}
Following the idea of trying to get the Fourier transforms of the
Riesz products in $\nu_\tau$ to cancel on $\{m.\Lambda: |m|=r\}$,
we should like a measure $\tau_l$ with $\|\tau_l\|$ minimal
subject to
\begin{equation*}
\int{t^kd\tau_l(t)}=0 \textrm{ for } 1<k\leq l,
\int{d\tau_l(t)}=0, \textrm{ and } \int{t d\tau_l(t)} = 1.
\end{equation*}
M\'{e}la, in \cite{JFM}, already had this idea, and moreover for
the purpose of constructing essentially the auxiliary measure we
want. To produce $\tau_l$ he constructs a measure $\sigma_l$ with
the following properties:
\begin{lemma}\emph{(Lemma 4, \S7, \cite{JFM})}
Suppose that $l>1$ is an integer. Then there is a measure
$\sigma_l$ on $[0,1]$ such that
\begin{equation*}
\int{s^{2k-1} d\sigma_l(s)}=0 \textrm{ for } 2 \leq k \leq l,
\int{s d\sigma_l(s)}=1 \textrm{ and } \|\sigma_l\| = 2l-1.
\end{equation*}
\end{lemma}
He chooses $\sigma_l$ to be (the measure induced by) the polynomial
$s+((-1)^l/(2l-1))P_{2l-1}(s)$ where $P_{2l-1}$ is the Chebychev
polynomial of order $2l-1$. Once this is known it is not hard to
verify the properties of $\sigma_l$.

We take $\tau_{2l}$ to be the odd measure on $[-1,1]$ which
extends $2\sigma_l(2s)$ on $[0,1/2]$, and the null measure on
$[1/2,1]$. It is easy, then, to verify the following.
\begin{lemma}\label{eqnlop}
Suppose that $l>1$ is an integer. Then the measure $\tau_{2l}$ on
$[-1,1]$ has $\|\tau_{2l}\| \leq 2(2l-1)$,
\begin{equation*}
\int{t^kd\tau_{2l}(t)} = \begin{cases} 0 \textrm{ if } k \leq 2l \textrm{ and } k \neq 1\\
1 \textrm{ if } k=1
\end{cases}
\end{equation*}
and $|\int{t^kd\tau_{2l}(t)}| \leq 2^{1-k}$ for all $k$.
\end{lemma}
Proposition \ref{Mela} follows from this by taking
$\mu_\eta=\nu_{\tau_{2l}}$ with $l=\lceil 2^{-1}\log_2
\eta^{-1}\rceil$.

\section{Remarks on Theorem
\ref{AGChang}}\label{remstruct}

The technique of applying Lemma \ref{maxdis} to reduce Theorem
\ref{AGChang} to Proposition \ref{dissprop} is used by Chang,
\cite{MCC}, in the proof of Theorem \ref{chang}. The analogue of
Proposition \ref{dissprop} in that case is proved using the dual
formulation of Rudin's inequality, which states that if $\Lambda$
is a dissociated set of characters on $G$ and $f \in L^2(G)$ then
\begin{equation*}
\|\wh{f}|_{\Lambda}\|_2 \ll \sqrt{\frac{p}{p-1}}\|f\|_p \textrm{
for } 2 \geq p>1.
\end{equation*}
Hal\'{a}sz, \cite{GH}, uses the inner product technique of
Proposition \ref{dissprop} to prove a non-Fourier result in
discrepancy theory and employs a Riesz product (for a different
Hilbert space) as the auxiliary measure. An exposition of his
result may be found in Chazelle \cite{BC} and this was the
original motivation for our result.

Green pointed out the fact that M\'{e}la, in \cite{JFM}, had
already used the method of linear combinations of Riesz products
to construct the auxiliary measure we require. M\'{e}la uses it as
an example to show that a result of his regarding
$\epsilon$-idempotent measures is essentially best possible. In
fact it follows from M\'{e}la's result that essentially no better
auxiliary measure than the one we have constructed exists.

Finally, Proposition \ref{dissprop}, like its analogue for Chang's
theorem, can be proved for a wider class of sets than simply
dissociated sets, namely Sidon sets. If $G$ is a compact Abelian
group and $\Lambda$ a set of characters on $G$ then $\Lambda$ is a
\emph{Sidon set} if every $\phi \in \ell^\infty(\Lambda)$ has the
form $\wh{\mu}|_\Lambda$ for some $\mu \in M(G)$. Extending the
results to Sidon sets comes down to constructing auxiliary
measures for this class of sets. Drury, in \cite{SD}, proved
Proposition \ref{primitiveauxiliary} for Sidon sets of characters
on general compact Abelian groups, and his methods can be used to
extend Proposition \ref{auxiliary}.

\section{Local Fourier analysis in compact Abelian
groups}\label{localft}

Our attention now turns to developing the tools of local Fourier
analysis and localizing the results of the previous four sections.

Given $f \in L^1(G)$ we often want to approximate $f$ by a less
complicated function. One way to do this is to approximate $f$ by
its expectation on level sets of characters. To analyze the error
in so doing we restrict the function to these level sets and use
the Fourier transform on the restricted function.

If $\Gamma$ is a finite set of characters then we define the
\emph{annihilator} of $\Gamma$ to be
\begin{equation*}
\Gamma^\perp:=\{x \in G: \gamma(x)=1 \textrm{ for all } \gamma \in
\Gamma\}.
\end{equation*}
If $x' + \Gamma^\perp$ (a maximal joint level set of the
characters in $\Gamma$) has positive measure in $G$ then it is
easy to localize the Fourier transform to $x'+\Gamma^{\perp}$:
\begin{equation*}
L^1(x'+\Gamma^\perp) \rightarrow \ell^\infty(\wh{G}); f \mapsto
\wh{fd(x'+\mu_{\Gamma^\perp})}.
\end{equation*}
Note that the right hand side is constant on cosets of
$\Gamma^{\perp\perp}$ (defined in the obvious manner) and so can
be thought of as an element of
$\ell^\infty(\wh{G}/\Gamma^{\perp\perp})$.

Bourgain, in \cite{JB}, observed that one can localize the Fourier
transform to typical approximate level sets and retain approximate
versions of a number of the standard results for the Fourier
transform on compact Abelian groups. Since his original work
various expositions and extensions have appeared most notably in
the various papers of Green and Tao. Indeed all the results of
this section can be found in \cite{BJGTCTU3}, for example.

\subsection{Approximate annihilators: Bohr neighborhoods and some of their properties}

Throughout this section $G$ is a compact Abelian group, $\Gamma$ a
non-empty finite set of characters on $G$ and $\delta \in (0,1]$.

We can define a natural valuation on $S^1$, namely
\begin{equation*}
\|z\|:=\frac{1}{2\pi}\inf_{n \in \mathbb{Z}}{| 2\pi n + \arg z|},
\end{equation*}
which can be used to measure how far $\gamma(x)$ is from 1.
Consequently we define a prototype for an approximate annihilator:
\begin{equation*}
B(\Gamma,\delta):=\{x \in G: \|\gamma(x)\| \leq \delta \textrm{
for all }\gamma \in \Gamma\},
\end{equation*}
called a \emph{Bohr set}. A translate of such a set is called a
\emph{Bohr neighborhood}. We adopt the convention that if
$B(\Gamma,\delta)$ is a Bohr set then the size of $\Gamma$ is
denoted by $d$.

Bohr sets are easily seen to be closed. To ensure that they have
positive measure we recall the following easy application of the
pigeonhole principle.
\begin{lemma}\label{bohrsize}
Suppose $G$ is a compact Abelian group and $B(\Gamma,\delta)$ is a
Bohr set. Then $\mu_G(B(\Gamma,\delta)) \geq \delta^d$, where as
our convention states, $d:=|\Gamma|$.
\end{lemma}
Hence we write $\beta_{\Gamma,\delta}$, or simply $\beta$ or
$\beta_\delta$ if the parameters are implicit, for the measure
induced on $B(\Gamma,\delta)$ by $\mu_G$, normalised so that
$\|\beta_{\Gamma,\delta}\|_1=1$. This is sometimes referred to as
the \emph{normalised Bohr cutoff}. We write $\beta'$ for
$\beta_{\Gamma',\delta'}$, or $\beta_{\Gamma,\delta'}$ if no
$\Gamma'$ has been defined.

Annihilators are subgroups of $G$, a property which, at least in
an approximate form, we would like to recover. Suppose that $\eta
\in (0,1]$. Then $B({\Gamma},\delta)+B({\Gamma},\eta\delta)
\subset B({\Gamma},(1+\eta)\delta)$. If
$B({\Gamma},(1+\eta)\delta)$ is not much bigger than
$B({\Gamma},\delta)$ then we have a sort of approximate additive
closure in the sense that
$B({\Gamma},\delta)+B({\Gamma},\eta\delta) \approx
B({\Gamma},(1+\eta)\delta)$. Not all Bohr sets have this property,
however, Bourgain showed that typically they do. For our purposes
we have the following proposition.
\begin{proposition}\label{ubreg}
Suppose that $G$ is a compact Abelian group, $\Gamma$ a set of $d$
characters on $G$ and $\delta \in (0,1]$. There is an absolute
constant $c_{\mathcal{R}}>0$ and a $\delta' \in [\delta/2,\delta)$
such that
\begin{equation*}
\frac{\mu_G(B(\Gamma,(1+\kappa)\delta'))}{\mu_G(B(\Gamma,\delta'))}
= 1 + O(|\kappa|d)
\end{equation*}
whenever $|\kappa|d \leq c_{\mathcal{R}}$.
\end{proposition}
This result is not as easy as the rest of the section, it uses a
covering argument; a nice proof can be found in \cite{BJGTCTU3}. We
say that $\delta'$ is \emph{regular for $\Gamma$} or that
$B(\Gamma,\delta')$ is \emph{regular} if
\begin{equation*}
\frac{\mu_G(B(\Gamma,(1+\kappa)\delta'))}{\mu_G(B(\Gamma,\delta'))}
= 1 + O(|\kappa|d) \textrm{ whenever } |\kappa|d \leq
c_{\mathcal{R}}.
\end{equation*}
It is regular Bohr sets to which we localize the Fourier transform
and we begin by observing that regular Bohr cutoffs are
approximately translation invariant and so function as normalised
approximate Haar measures.
\begin{lemma}
\label{contlem}\emph{(Normalized approximate Haar measure)}
Suppose that $G$ is a compact Abelian group and $B(\Gamma,\delta)$
is a regular Bohr set. If $y \in B(\Gamma,\delta')$ then
$\|(y+\beta_\delta) - \beta_\delta\| \ll d\delta'\delta^{-1}$
where we recall that $y+\beta_\delta$ denotes the measure
$\beta_\delta$ composed with translation by $y$.
\end{lemma}
\begin{proof}
Note that $\supp{((y+\beta_{\delta}) - \beta_{\delta})} \subset
B(\Gamma,\delta+\delta') \setminus B(\Gamma,\delta-\delta')$ whence
\begin{equation*}
\|(y+\beta_{\delta}) - \beta_{\delta}\| \leq
\frac{\mu_G(B(\Gamma,\delta+\delta') \setminus
B(\Gamma,\delta-\delta'))}{\mu_G(B(\Gamma,\delta))} \ll
d\delta'\delta^{-1},
\end{equation*}
by regularity.
\end{proof}
In applications the following two simple corollaries will be
useful but they should be ignored until they are used.
\begin{corollary}
\label{contlemcor} Suppose that $G$ is a compact Abelian group and
$B(\Gamma,\delta)$ is a regular Bohr set. If $\mu \in
M(B(\Gamma,\delta'))$ then $\|\beta\ast \mu - \beta\int{d\mu}\|
\ll \|\mu\|d\delta'\delta^{-1}$.
\end{corollary}
\begin{proof} The measures $\beta \ast \mu$ and $\beta \int{d\mu}$ agree inside
$B(\Gamma,\delta-\delta')$ and outside $B(\Gamma,\delta+\delta')$,
furthermore $\|\beta \ast \mu\| \leq \|\mu\|$ and $\|\beta
\int{d\mu}\| \leq \|\mu\|$, whence
\begin{equation*}
\|\beta \ast \mu - \beta \int{d\mu}\| \ll
\|\mu\|\frac{\mu_G(B(\Gamma,\delta+\delta') \setminus
B(\Gamma,\delta - \delta'))}{\mu_G(B(\Gamma,\delta))} \ll
\|\mu\|d\delta'\delta^{-1}
\end{equation*}
by regularity.
\end{proof}
\begin{corollary}
\label{contlemcor2} Suppose that $G$ is a compact Abelian group
and $B(\Gamma,\delta)$ is a regular Bohr set. If $f \in
L^\infty(G)$ then
\begin{equation*}
\|f\ast \beta - f \ast \beta(x)\|_{L^\infty(x+B(\Gamma,\delta'))}
\ll \|f\|_{L^\infty(G)}d\delta'\delta^{-1}.
\end{equation*}
\end{corollary}
\begin{proof}
Note that
\begin{eqnarray*}
|f \ast \beta(x+y) - f \ast \beta(x)| & = & |f \ast
((-y+\beta)-\beta)(x)|\\ & \leq &
\|f\|_{L^\infty(\mu_G)}\|(-y+\beta)-\beta\|.
\end{eqnarray*}
The result follows by Lemma \ref{contlem}.
\end{proof}

With an approximate Haar measure we are in a position to define
the local Fourier transform: Suppose that $x'+B(\Gamma,\delta)$ is
a regular Bohr neighborhood. Then we define the Fourier transform
local to $x'+B(\Gamma,\delta)$ by
\begin{equation*}
L^1(x'+B(\Gamma,\delta)) \rightarrow \ell^\infty(\wh{G}); f
\mapsto \wh{fd(x'+\beta_{\Gamma,\delta})},
\end{equation*}
where $L^p(x'+B(\Gamma,\delta))$ denotes the space
\begin{equation*}
\{f \in L^1(G):\supp f \subset x'+B(\Gamma,\delta) \textrm{ and }
\int{|f|^pd(x'+\beta)} < \infty \},
\end{equation*}
equipped with the norm
\begin{equation*}
\|f\|_{L^p(x'+B(\Gamma,\delta))}:=\left(\int{|f|^pd(x'+\beta)}\right)^{1/p}.
\end{equation*}
It is useful at this stage to also define $A(x'+B(\Gamma,\delta))$
which is the space
\begin{equation*}
\{f \in L^1(G):\supp f \subset x'+B(\Gamma,\delta) \textrm{ and }
\|\wh{f}\|_1 < \infty \},
\end{equation*}
equipped with the norm
\begin{equation*}
\|f\|_{A(x'+B(\Gamma,\delta))}=\|\wh{f}\|_1.
\end{equation*}
The translation of the Bohr set by $x'$ simply twists the Fourier
transform and is unimportant for the most part so we tend to
restrict ourselves to the case when $x'=0$.

$\wh{fd\mu_{\Gamma^\perp}}$ was constant on cosets of
$\Gamma^{\perp\perp}$. In the approximate setting have an
approximate analogue on which $\wh{gd\beta}$ does not vary too
much. There are a number of possibilities:
\begin{eqnarray*}
&\{\gamma:|1-\gamma(x)| \leq \epsilon \textrm{ for all
} x \in B(\Gamma,\delta)\} &\textrm{ for } \epsilon \in (0,1] \\
&\{\gamma:|1-\wh{\beta}(\gamma)| \leq \epsilon\}
&\textrm{ for } \epsilon  \in (0,1]\\
&\{\gamma:|\wh{\beta}(\gamma)| \geq \epsilon \} &\textrm{ for }
\epsilon  \in (0,1].
\end{eqnarray*}
In applications each of these classes of sets is useful and so we
should like all of them to be approximately equivalent. There is a
clear chain of inclusions between the classes:
\begin{equation*}
\{\gamma:|1-\gamma(x)| \leq \epsilon \textrm{ for all } x \in
B(\Gamma,\delta)\} \subset \{\gamma:|1-\wh{\beta}(\gamma)| \leq
\epsilon\} \subset \{\gamma:|\wh{\beta}(\gamma)| \geq 1-\epsilon
\}
\end{equation*}
for $\epsilon \in (0,1]$. For a small cost in the width of the
Bohr set we can ensure that the sets in the third class are
contained in those in the first.
\begin{lemma} \label{nestsupport}
Suppose that $G$ is a compact Abelian group and $B(\Gamma,\delta)$
is a regular Bohr set. Suppose that $\eta_1,\eta_2>0$. Then there
is a $\delta' \gg \eta_1\eta_2\delta/d$ such that
\begin{equation*}
\{\gamma: |\wh{\beta}(\gamma)| \geq \eta_1\} \subset \{\gamma:
|1-\gamma(x)| \leq \eta_2 \textrm{ for all } x \in
B(\Gamma,\delta')\}.
\end{equation*}
\end{lemma}
\begin{proof}
If $x \in B(\Gamma,\delta')$ then we have
\begin{equation*}
\eta_1|1-\gamma(x)| \leq
|\wh{\beta}(\gamma)||1-\overline{\gamma(x)}|=|\wh{((x+\beta) -
\beta)}(\gamma)| \ll d \delta'\delta^{-1}
\end{equation*}
by Lemma \ref{contlem}. It follows that we may pick $\delta' \gg
\eta_1\eta_2\delta/d$ such that $|1-\gamma(x)| \leq \eta_2$ for all
$x \in B(\Gamma,\delta')$.
\end{proof}

This concludes the basic definitions of the local Fourier
transform; in the next section we transfer the results we require
to the local setting.

\section{A structural result for the local Fourier spectrum}\label{localstruct}

In \S\ref{struct} we recorded a number of results regarding the
structure of the collection of characters supporting large values
of the Fourier transform; in this section we examine local
versions of these. About the simplest statement we made was
(\ref{bessel}) which asserted that if $f \in L^2(G)$ and $\Gamma =
\{\gamma: |\wh{f}(\gamma)| \geq \epsilon\|f\|_1\}$ then
\begin{equation*}
|\Gamma| \leq \epsilon^{-2} (\|f\|_2\|f\|_1^{-1})^2.
\end{equation*}
The following analogue for functions $f \in L^2(B(\Gamma,\delta))$
was proved by Green and Tao in \cite{BJGTCTU3}, by localizing
Bessel's inequality.
\begin{proposition}\label{locbes}
Suppose that $G$ is a compact Abelian group and $B(\Gamma,\delta)$ a
regular Bohr set. Suppose that $f \in L^2(B(\Gamma,\delta))$ and
$\epsilon,\eta \in (0,1]$. Write $L_f$ for the quantity
$\|f\|_{L^2(B(\Gamma,\delta))}\|f\|_{L^1(B(\Gamma,\delta))}^{-1}$.
Then there is a set $\Lambda$ of characters and a $\delta' \in
(0,1]$ with
\begin{equation*}
|\Lambda| \ll \epsilon^{-2}L_f^2 \textrm{ and } \delta' \gg
\epsilon^2\eta\delta/dL_f^2,
\end{equation*}
such that
\begin{equation*}
\{\gamma \in \wh{G}:|\wh{fd\beta}(\gamma)| \geq \epsilon
\|f\|_{L^1(B(\Gamma,\delta))}\} \subset \{\gamma: |1-\gamma(x)|
\leq \eta \textrm{ for all } x \in B(\Gamma \cup
\Lambda,\delta')\}.
\end{equation*}
\end{proposition}
Later on in \S\ref{struct} we noted (\ref{l1bessel}), an analogue
of (\ref{bessel}) for functions $f \in A(G)$, and there is a
corresponding local analogue of this result for $f \in
A(B(\Gamma,\delta))$ which we note now.
\begin{proposition}\label{localag}
Suppose that $G$ is a compact Abelian group and $B(\Gamma,\delta)$ a
regular Bohr set. Suppose that $f \in A(B(\Gamma,\delta))$ and
$\epsilon,\eta \in (0,1]$. Write $A_f$ for the quantity
$\|f\|_{A(B(\Gamma,\delta))}\|f\|_{L^\infty(B(\Gamma,\delta))}^{-1}$.
Then there is a set $\Lambda$ of characters and a $\delta' \in
(0,1]$ with
\begin{equation*}
|\Lambda| \ll \epsilon^{-1}A_f \textrm{ and } \delta' \gg
\epsilon\eta\delta/dA_f,
\end{equation*}
such that
\begin{equation*}
\{\gamma \in \wh{G}:|\wh{fd\beta}(\gamma)| \geq \epsilon
\|f\|_{L^\infty(B(\Gamma,\delta))}\} \subset \{\gamma:
|1-\gamma(x)| \leq \eta \textrm{ for all } x \in B(\Gamma \cup
\Lambda,\delta')\}.
\end{equation*}
\end{proposition}
In \cite{TSASS} the result of Green and Tao is refined with
Chang's theorem. One can use the same techniques to refine
Proposition \ref{localag} with Theorem \ref{AGChang}, our
$A(G)$-analogue of Chang's theorem; doing so gives the following.
\begin{proposition}\label{f2agchang}
Suppose that $G$ is a compact Abelian group and $B(\Gamma,\delta)$ a
regular Bohr set. Suppose that $f \in A(B(\Gamma,\delta))$ and
$\epsilon,\eta \in (0,1]$. Write $A_f$ for the quantity
$\|f\|_{A(B(\Gamma,\delta))}\|f\|_{L^\infty(B(\Gamma,\delta))}^{-1}$.
Then there is a set $\Lambda$ of characters and a $\delta' \in
(0,1]$ with
\begin{equation*}
|\Lambda| \ll \epsilon^{-1}(1+\log A_f) \textrm{ and } \delta' \gg
\epsilon\eta\delta/d^2(1+\log A_f),
\end{equation*}
such that
\begin{equation*}
\{\gamma \in \wh{G}:|\wh{fd\beta}(\gamma)| \geq \epsilon
\|f\|_{L^\infty(B(\Gamma,\delta))}\} \subset \{\gamma:
|1-\gamma(x)| \leq \eta \textrm{ for all } x \in B(\Gamma \cup
\Lambda,\delta')\}.
\end{equation*}
\end{proposition}
For our application we are in fact able to assume that $f \in A(G)$,
and in that case we have the following slightly more general
conclusion.
\begin{proposition}\label{localagchang}
Suppose that $G$ is a compact Abelian group and $B(\Gamma,\delta)$ a
regular Bohr set. Suppose that $f \in A(G)$ and $\epsilon,\eta \in
(0,1]$. Write $A_f$ for the quantity
$\|f\|_{A(G)}\|f\|_{L^\infty(B(\Gamma,\delta))}^{-1}$. Then there is
a set $\Lambda$ of characters and a $\delta' \in (0,1]$ with
\begin{equation*}
|\Lambda| \ll \epsilon^{-1}(1+\log A_f) \textrm{ and } \delta' \gg
\epsilon^2\eta\delta/d^2(1+\log A_f),
\end{equation*}
such that
\begin{equation*}
\{\gamma \in \wh{G}:|\wh{fd\beta}(\gamma)| \geq \epsilon
\|f\|_{L^\infty(B(\Gamma,\delta))}\} \subset \{\gamma:
|1-\gamma(x)| \leq \eta \textrm{ for all } x \in B(\Gamma \cup
\Lambda,\delta')\}.
\end{equation*}
\end{proposition}
We do not require Propositions \ref{localag} or \ref{f2agchang} in
our application and their methods of proof are in both cases
simplifications of the method used for Proposition
\ref{localagchang} so we shall restrict our attention to proving
that result.

A key tool in Theorem \ref{AGChang} is that of dissociativity; in
the local setting we use the following version of it. If $S$ is a
non-empty symmetric neighborhood of $0_{\wh{G}}$ then we say that
$\Lambda$ is \emph{$S$-dissociated} if
\begin{equation*}
m \in \{-1,0,1\}^\Lambda \textrm{ and } m.\Lambda \in S \textrm{
implies that } m\equiv 0.
\end{equation*}
Vanilla dissociativity corresponds to taking $S=\{0_{\wh{G}}\}$,
and typically in the local setting $S$ will be a set of characters
at which $\wh{\beta}$ is large for some Bohr set
$B(\Gamma,\delta)$.

In the same way as Theorem \ref{AGChang} follows from Lemma
\ref{maxdis} and Proposition \ref{dissprop}, Proposition
\ref{localagchang} follows from the next two lemmas.
\begin{lemma}\label{spanning agchang}
Suppose that $G$ is a compact Abelian group and $B(\Gamma,\delta)$
is a regular Bohr set. Suppose that $\eta',\eta \in (0,1]$ and
$\Delta$ is a set of characters on $G$. If $\Lambda$ is a maximal
$\{\gamma:|\wh{\beta}(\gamma)| \geq \eta' \}$-dissociated subset
of $\Delta$ then there is a $\delta' \gg
\min\{\eta/|\Lambda|,\eta' \eta \delta/d\}$ such that
\begin{equation*}
\Delta \subset \{\gamma:|1-\gamma(x)| \leq \eta \textrm{ for all }
x \in B(\Gamma \cup \Lambda,\delta')\}.
\end{equation*}
\end{lemma}
\begin{lemma}\label{content agchang}
Suppose that $G$ is a compact Abelian group and $B(\Gamma,\delta)$
is a regular Bohr set. Suppose that $f \in A(G)$ and $\epsilon,\eta
\in (0,1]$. Write $A_f$ for the quantity
$\|f\|_{A(G)}\|f\|_{L^\infty(B(\Gamma,\delta))}^{-1}$. Then there is
a $\delta' \gg \epsilon^2\delta/d(1+\log A_f)$ regular for $\Gamma$
such that if $\Lambda$ is a $\{\gamma:|\wh{\beta'}(\gamma)| \geq
1/3\}$-dissociated subset of $\{\gamma \in
\wh{G}:|\wh{fd\beta}(\gamma)| \geq \epsilon
\|f\|_{L^\infty(B(\Gamma,\delta))}\}$ then
\begin{equation*}
|\Lambda| \ll \epsilon^{-1} (1+\log A_f).
\end{equation*}
\end{lemma}

\subsection{The proof of Lemma \ref{spanning agchang}}

The following lemma localizes Lemma \ref{maxdis}, which
corresponds to the case $S=\{0_{\wh{G}}\}$.
\begin{lemma}
Suppose that $G$ is a compact Abelian group. Suppose that $S$ is a
non-empty symmetric neighborhood of $0_{\wh{G}}$. Suppose that
$\Delta$ is a set of characters on $G$ and $\Lambda$ is a maximal
$S$-dissociated subset of $\Delta$. Then $\Delta \subset \langle
\Lambda \rangle + S$.
\end{lemma}
\begin{proof}
If $\lambda_0 \in \Delta \setminus (\langle \Lambda \rangle +S)$
then we put $\Lambda':=\Lambda \cup \{ \lambda_0 \}$, which is a
strict superset of $\Lambda$, and a subset of $\Delta$. It turns
out that $\Lambda'$ is also $S$-dissociated which contradicts the
maximality of $\Lambda$. Suppose that $m:\Lambda' \rightarrow
\{-1,0,1\}$ and $m.\Lambda' \in S$. Then we have three
possibilities for the value of $m_{\lambda_0}$:
\begin{enumerate}
\item $m.\Lambda' = \lambda_0 + m|_\Lambda.\Lambda$, in which case
$\lambda_0 \in -m|_\Lambda.\Lambda + S \subset \langle \Lambda
\rangle + S$ - a contradiction; \item $m.\Lambda' = - \lambda_0 +
m|_\Lambda.\Lambda$, in which case $\lambda_0 \in
m|_\Lambda.\Lambda - S \subset \langle \Lambda \rangle + S$ - a
contradiction; \item $m.\Lambda' = m|_\Lambda.\Lambda$, in which
case $m|_\Lambda\equiv 0$ since $\Lambda$ is $S$-dissociated and
hence $m \equiv 0$.
\end{enumerate}
It follows that $m.\Lambda' \in S\Rightarrow m\equiv 0$ i.e.
$\Lambda'$ is $S$-dissociated as claimed.
\end{proof}
Lemma \ref{spanning agchang} then follows from the above and the
next lemma.
\begin{lemma}
Suppose that $G$ is a compact Abelian group and $B(\Gamma,\delta)$
is a regular Bohr set. Suppose that $\eta',\eta \in (0,1]$ and
$\Lambda$ is a set of characters on $G$. Then there is a $\delta'
\gg \min\{\eta/|\Lambda|,\eta'\eta\delta/d\}$ such that
\begin{equation*}
\langle \Lambda \rangle + \{\gamma:|\wh{\beta}(\gamma)| \geq \eta'
\} \subset \{\gamma:|1-\gamma(x)| \leq \eta \textrm{ for all } x
\in B(\Gamma \cup \Lambda,\delta')\}.
\end{equation*}
\end{lemma}
\begin{proof}
The lemma has two parts.
\begin{enumerate}
\item If $\lambda \in \langle \Lambda \rangle$ then
\begin{equation*}
|1-\lambda(x)| \leq \sum_{\lambda' \in \Lambda}{|1-\lambda'(x)|},
\end{equation*}
so there is a $\delta'' \gg \eta/|\Lambda|$ such that
\begin{equation*}
\langle \Lambda \rangle \subset \{\gamma: |1-\gamma(x)| \leq
\eta/2 \textrm{ for all } x \in B(\Lambda,\delta'')\}.
\end{equation*}
\item By Lemma \ref{nestsupport} there is a $\delta''' \gg
\eta\eta'\delta/d$ such that
\begin{equation*}
\{\gamma:|\wh{\beta}(\gamma)| \geq \eta' \} \subset
\{\gamma:|1-\gamma(x)| \leq \eta/2 \textrm{ for all } x \in
B(\Gamma,\delta''')\}.
\end{equation*}
\end{enumerate}
Taking $\delta' = \min\{\delta'',\delta'''\}$ we have the result
by the triangle inequality.
\end{proof}

\subsection{The proof of Lemma \ref{content agchang}}
The proof follows that of Proposition \ref{dissprop} with the
additional ingredient of smoothed measures.

Suppose that $B(\Gamma,\delta)$ is a regular Bohr set on $G$. For
$L \in \mathbb{N}$ and $\kappa \in (0,1]$ we write
\begin{equation*}
\tilde{\beta}^{L,\kappa}_{\Gamma,\delta}:=\beta_{\Gamma,(1-\kappa)\delta}
\ast \beta_{\Gamma,\kappa\delta/L}^L,
\end{equation*}
where here, as in future, juxtaposition of measures denotes
convolution. $\tilde{\beta}^{L,\kappa}_{\Gamma,\delta}$ is a good
approximation to $\beta$ in $M(G)$:
\begin{equation*}
\|\tilde{\beta}_{\Gamma,\delta}^{L,\kappa}-\beta_{\Gamma,\delta}\|
\leq \|\beta_{\Gamma,\delta(1-\kappa)} \ast \mu -
\beta_{\Gamma,\delta(1-\kappa)}\| + \|
\beta_{\Gamma,\delta(1-\kappa)} -\beta_{\Gamma,\delta}\|
\end{equation*}
where $\mu =\beta_{\Gamma,\kappa\delta/L}^L$, the convolution of
$\beta_{\Gamma,\kappa\delta/L}$ with itself $L$ times. We deal
with the first term using Corollary \ref{contlemcor} which yields
\begin{equation*}
\|\beta_{\Gamma,\delta(1-\kappa)} \ast \mu -
\beta_{\Gamma,\delta(1-\kappa)}\| \ll \kappa d,
\end{equation*}
since $\supp \mu \subset B(\Gamma,\kappa\delta)$. For the second
term we have
\begin{eqnarray*}
\| \beta_{\Gamma,\delta(1-\kappa)} -\beta_{\Gamma,\delta}\| & \leq
& \| \beta_{\Gamma,\delta(1-\kappa)}
-\beta_{\Gamma,\delta}|_{B(\Gamma,\delta(1-\kappa))}\| +
\|\beta_{\Gamma,\delta}|_{B(\Gamma,\delta)\setminus
B(\Gamma,\delta(1-\kappa))}\|\\ & = & \left(1-
\frac{\mu_G(B(\Gamma,\delta(1-\kappa)))}{\mu_G(B(\Gamma,\delta))}\right)
+ \|\beta_{\Gamma,\delta}|_{B(\Gamma,\delta)\setminus
B(\Gamma,\delta(1-\kappa))}\|\\ & = & O(\kappa d) +
\|\beta_{\Gamma,\delta}|_{B(\Gamma,\delta)\setminus
B(\Gamma,\delta(1-\kappa))}\| \textrm{ by regularity,}\\ & = &
O(\kappa d) + \left(\frac{\mu_G(B(\Gamma,\delta)) -
\mu_G(B(\Gamma,\delta(1-\kappa)))}{\mu_G(B(\Gamma,\delta))}\right)\\
& = & O(\kappa d) \textrm{ by regularity.}
\end{eqnarray*}
It follows that $\|\tilde{\beta}_{\Gamma,\delta}^{L,\kappa} -
\beta\|=O(\kappa d)$ and hence if $f \in
L^\infty(B(\Gamma,\delta))$ we have
\begin{equation}\label{smoothdi}
|\wh{fd\tilde{\beta}^{L,\kappa}_{\Gamma,\delta}}(\gamma) -
\wh{fd\beta_{\Gamma,\delta}}(\gamma)| \ll
\|f\|_{L^\infty(B(\Gamma,\delta))}\kappa d.
\end{equation}

\begin{proof} We begin by fixing $\kappa$
and $L$ in the smoothed measure
$\tilde{\beta}^{L,\kappa}_{\Gamma,\delta}$ so that we may dispense
with the superscripts and subscripts and simply write
$\tilde{\beta}$. Take $L=2R$, where $R$ will be chosen later and
$\kappa \gg \epsilon/d$ so that
\begin{equation*}
|\wh{fd\tilde{\beta}^{L,\kappa}_{\Gamma,\delta}}(\gamma) -
\wh{fd\beta_{\Gamma,\delta}}(\gamma)| \leq
2^{-1}\epsilon\|f\|_{L^\infty(B(\Gamma,\delta))},
\end{equation*}
which we may certainly do by (\ref{smoothdi}), and also so that
$\delta':=\kappa\delta/L$ is regular for $\Gamma$. As usual this
last requirement is possible by Proposition \ref{ubreg}. It
follows that
\begin{equation*}
|\wh{fd\beta}(\gamma)| \geq \epsilon
\|f\|_{L^\infty(B(\Gamma,\delta))} \Rightarrow
|\wh{fd\tilde{\beta}}(\gamma)| \geq
2^{-1}\epsilon\|f\|_{L^\infty(B(\Gamma,\delta))}.
\end{equation*}
Now suppose that $\Lambda$ is a $\{\gamma:|\wh{\beta'}(\gamma)| \geq
1/3\}$-dissociated subset of $\{\gamma \in \wh{G}:
|\wh{fd\beta}\gamma)| \geq \epsilon
\|f\|_{L^\infty(B(\Gamma,\delta))}\}$ and that $\Lambda' \subset
\Lambda$ has size at most $R$. $\Lambda'$ is certainly still
$\{\gamma:|\wh{\beta'}(\gamma)| \geq 1/3\}$-dissociated. We define
\begin{equation*}
\omega(\lambda):=\frac{\wh{fd\tilde{\beta}}(\lambda)}{|\wh{fd\tilde{\beta}}(\lambda)|}
\textrm{ for all }\lambda \in \Lambda'.
\end{equation*}
$\omega \in \ell^\infty(\Lambda')$, $\|\omega\|_\infty \leq 1$,
$\Lambda'$ is finite and $\Lambda'$ is dissociated (since it is
$\{\gamma:|\wh{\beta'}(\gamma)| \geq 1/3\}$-dissociated) so we may
apply Proposition \ref{auxiliary} to get the auxiliary measure
$\mu_\eta$. To leverage the stronger dissociativity condition we
introduce a Riesz product:
\begin{equation*}
q:=\prod_{\lambda \in \Lambda'}{\left(1+\frac{\lambda +
\overline{\lambda}}{2}\right)}.
\end{equation*}
Recall (from \S\ref{RPs2} if necessary) that $q$ is non-negative and
since $\Lambda'$ is dissociated $\|q\|_1 =1$.

Plancherel's theorem gives
\begin{equation*}
\langle fd\tilde{\beta},\mu_\eta\ast q \rangle = \sum_{\gamma \in
\wh{G}}{\wh{fd\tilde{\beta}}(\gamma)\overline{\wh{q}(\gamma)\wh{\mu_\eta}(\gamma)}}.
\end{equation*}
We begin by bounding the right hand side from below using the bound
on $|\wh{\mu_\eta}(\lambda)|$ for $\lambda \not\in \wh{G}$ and the
fact that $\wh{q}(\lambda) \geq 1/2$ if $\lambda \in \Lambda'$.
\begin{eqnarray*}
|\sum_{\gamma \in
\wh{G}}{\wh{fd\tilde{\beta}}(\gamma)\overline{\wh{q}(\gamma)\wh{\mu_\eta}(\gamma)}}|
& = & |\sum_{\lambda \in
\Lambda'}{\wh{fd\tilde{\beta}}(\lambda)\overline{\wh{q}(\lambda)\wh{\mu_\eta}(\lambda)}}+
\sum_{\lambda \not \in
\Lambda'}{\wh{fd\tilde{\beta}}(\lambda)\overline{\wh{q}(\lambda)\wh{\mu_\eta}(\lambda)}}|\\
& \geq & |\sum_{\lambda \in
\Lambda'}{\wh{fd\tilde{\beta}}(\lambda)\overline{\wh{q}(\lambda)\wh{\mu_\eta}(\lambda)}}|
- |\sum_{\lambda \not \in
\Lambda'}{\wh{fd\tilde{\beta}}(\lambda)\overline{\wh{q}(\lambda)\wh{\mu_\eta}(\lambda)}}|\\
& \geq & |\sum_{\lambda \in
\Lambda'}{\wh{fd\tilde{\beta}}(\lambda)\overline{\wh{q}(\lambda)\omega(\lambda)}}|
- \eta\sum_{\lambda \in
\wh{G} }{|\wh{q}(\lambda)||\wh{fd\tilde{\beta}}(\lambda)|}\\
& \geq & 2^{-1}\sum_{\lambda \in
\Lambda'}{|\wh{fd\tilde{\beta}}(\lambda)|} - \eta\sum_{\lambda \in
\wh{G}}{\sum_{\gamma
\in\wh{G}}{|\wh{q}(\lambda)||\wh{f}(\gamma)\wh{\tilde{\beta}}(\lambda - \gamma)|}}\\
& \geq & 2^{-1}\sum_{\lambda \in
\Lambda'}{|\wh{fd\tilde{\beta}}(\lambda)|} -
\eta\|f\|_{A(G)}\sup_{\gamma \in \wh{G}}{\sum_{\lambda \in
\wh{G}}{|\wh{q}(\lambda)||\wh{\tilde{\beta}}(\lambda - \gamma)|}}.
\end{eqnarray*}
For any $\gamma \in \wh{G}$ we can estimate the last sum in a
manner independent of $\gamma$ by using a positivity argument:
\begin{eqnarray*}
\sum_{\lambda \in
\wh{G}}{|\wh{q}(\lambda)||\wh{\tilde{\beta}}(\lambda - \gamma)|} &
= & \sum_{\lambda \in
\wh{G}}{|\wh{q}(\gamma-\lambda)||\wh{\tilde{\beta}}(\lambda)|}
\textrm{ by symmetry of $\wh{\tilde{\beta}}$,}\\ & = &
\sum_{\lambda \in
\wh{G}}{|\wh{q}(\gamma-\lambda)||\wh{\beta}(\lambda)
\wh{\beta'}(\lambda)^L|} \textrm{ by definition of
$\wh{\tilde{\beta}}$,}\\ & \leq & \sum_{\lambda \in \wh{G}
}{|\wh{q}(\gamma-\lambda)||\wh{\beta'}(\lambda)|^L} \textrm{ since
$|\wh{\beta}(\lambda)| \leq \|\beta\| = 1$, }\\ & = &
\wh{qd\beta'^L}(\gamma) \textrm{ since $L$ is even and $\wh{q}
\geq 0$,}\\ & \leq & \|qd\beta'^L\| \\ & =
&\wh{qd\beta'^L}(0_{\wh{G}}) \textrm{ by non-negativity of
$qd\beta'^L$,}
\\  & = & \sum_{\lambda \in
\wh{G}}{\wh{q}(\lambda)|\wh{\beta'}(\lambda)|^L} \textrm{ by
symmetry of $\wh{q}$.}
\end{eqnarray*}
We estimate this by splitting the range of summation into two
parts:
\begin{equation}\label{crty} \sum_{\lambda \in
\wh{G}}{\wh{q}(\lambda)|\wh{\beta'}(\lambda)|^L} \leq
\sum_{\lambda:|\wh{\beta'}(\lambda)| \geq
1/3}{\wh{q}(\lambda)|\wh{\beta'}(\lambda)|^L} +
\sum_{\lambda:|\wh{\beta'}(\lambda)| \leq
1/3}{\wh{q}(\lambda)|\wh{\beta'}(\lambda)|^L}.
\end{equation}
\begin{enumerate}
\item For the first sum: $|\wh{q}(\lambda)| \leq \|q\|_1 = 1$ and
$|\wh{\beta'}(\lambda)^L| \leq \|\beta'^L\| = 1$ so that each
summand is at most 1, furthermore $\supp \wh{q} \subset \langle
\Lambda' \rangle$ so
\begin{equation*}
\sum_{\lambda:|\wh{\beta'}(\lambda)| \geq
1/3}{\wh{q}(\lambda)|\wh{\beta'}(\lambda)|^L} \leq \sum_{\lambda \in
\langle \Lambda' \rangle:|\wh{\beta'}(\lambda)| \geq 1/3}{1}.
\end{equation*}
This range of summation contains at most 1 element by
$\{\gamma:|\wh{\beta'}(\gamma)|\geq 1/3\}$-dissociativity of
$\Lambda'$, and hence the sum is bounded above by 1. \item For the
second sum: $|\wh{q}(\lambda)| \leq \|q\|_1 =1$ and
$|\wh{\beta'}(\lambda)^L| \leq 3^{-L}$ for $\lambda$ in the range of
summation so that each summand is at most $9^{-|\Lambda'|}$, however
$\supp \wh{q} \subset \langle \Lambda' \rangle$ and $|\langle
\Lambda' \rangle| \leq 3^{|\Lambda'|}$ so
\begin{equation*}
\sum_{\lambda:|\wh{\beta'}(\lambda)| \leq
1/3}{\wh{q}(\lambda)|\wh{\beta'}(\lambda)|^L} \leq \sum_{\lambda \in
\langle \Lambda' \rangle}{9^{-|\Lambda'|}} \leq 1.
\end{equation*}
\end{enumerate}
It follows that the right hand side of (\ref{crty}) is bounded
above by 2, and working backwards these estimates combine to show
that
\begin{equation*}
\sum_{\lambda \in
\wh{G}}{|\wh{q}(\lambda)||\wh{\tilde{\beta}}(\lambda - \gamma)|}
\leq 2 \textrm{ for all } \gamma \in \wh{G},
\end{equation*}
and hence that
\begin{equation}\label{ytu}
|\langle fd\tilde{\beta},\mu_\eta\ast q \rangle| \geq
2^{-1}\sum_{\lambda \in \Lambda'}{|\wh{fd\tilde{\beta}}(\lambda)|} -
2\eta\|f\|_{A(G)}.
\end{equation}
To estimate the inner product from above we have:
\begin{equation*}
|\langle fd\tilde{\beta},\mu_\eta \ast q \rangle| \leq
\|f\|_{L^\infty(B(\Gamma,\delta))}\|\tilde{\beta}\|
\|\mu_\eta\|\|q\|_1 \ll \|f\|_{L^\infty(B(\Gamma,\delta))}(1+ \log
\eta^{-1})
\end{equation*}
by the estimate for $\|\mu_\eta\|$ given in Proposition
\ref{auxiliary}. Combining this with our lower bound for the inner
product in (\ref{ytu}) and the fact that if $\lambda \in \Lambda'$
then $|\wh{fd\tilde{\beta}}(\lambda)| \geq 2^{-1}\epsilon
\|f\|_{L^\infty(B(\Gamma,\delta))}$ gives
\begin{equation*}
\|f\|_{L^\infty(B(\Gamma,\delta))}(1+ \log  \eta^{-1}) +
\eta\|f\|_{A(G)} \gg |\Lambda'|\epsilon
\|f\|_{L^\infty(B(\Gamma,\delta))}.
\end{equation*}
Choosing
$\eta^{-1}=\|f\|_{A(G)}\|f\|_{L^\infty(B(\Gamma,\delta))}^{-1}$
yields that
\begin{equation*}
|\Lambda'| \ll \epsilon^{-1}(1+\log A_f).
\end{equation*}
Let $C$ be the absolute constant implicit in the notation on the
right so that $|\Lambda'| \leq C\epsilon^{-1}(1+\log A_f)$ is always
true, and set $R:=\lceil C\epsilon^{-1}(1+\log A_f) \rceil+1$. If
$|\Lambda|$ is a $\{\gamma:|\wh{\beta'}(\gamma)|\geq
1/3\}$-dissociated set of size greater than $R$, then let $\Lambda'$
be a subset of $\Lambda$ of size $R$, which is automatically
$\{\gamma:|\wh{\beta'}(\gamma)|\geq 1/3\}$-dissociated because
$\Lambda$ is $\{\gamma:|\wh{\beta'}(\gamma)|\geq 1/3\}$-dissociated.
By the above
\begin{equation*}
C\epsilon^{-1}(1+\log A_f)<\lceil C\epsilon^{-1}(1+\log A_f)
\rceil+1 = R=|\Lambda'| \leq C\epsilon^{-1}(1+\log A_f),
\end{equation*}
which is a contradiction and hence if $\Lambda$ is
$\{\gamma:|\wh{\beta'}(\gamma)|\geq 1/3\}$-dissociated then
$|\Lambda| < R \ll \epsilon^{-1}(1+\log A_f)$ as required.
\end{proof}

\section{An introduction to Littlewood's problem}\label{tfr}

Finally we turn to addressing the problem announced in the
introduction, but although the contents of \S\ref{localstruct} may
be taken as a black box for the purposes of the following
sections, \S\ref{localft} is a necessary notational prerequisite.

It is natural to begin by considering a qualitative analogue of
our problem, in particular we shall start by proving the following
well known result.
\begin{proposition}\label{qualitative Littlewood}\emph{(Qualitative Littlewood problem)}
Suppose that $A \subset \mathbb{T}$ has density $\alpha$ with $0 <
\alpha < 1$. Then $\chi_A \not \in A(\mathbb{T})$.
\end{proposition}
The proof of this proceeds in three stages, the first two of which
are naturally set in an arbitrary compact Abelian group $G$.
\begin{enumerate}
\item \emph{(Fourier inversion)} First, if $f \in A(G)$ then we
may define the function
\begin{equation*}
\wt{f}(x):=\sum_{\gamma \in \wh{G}}{\wh{f}(\gamma)\gamma(x)},
\end{equation*}
which is continuous since it is the uniform limit of continuous
functions. The Fourier inversion theorem tells us that
$\|\wt{f}-f\|_{L^\infty(G)}=0$. \item \emph{(Averaging)} Secondly,
by averaging there are elements $x_0,x_1 \in G$ such that
\begin{equation*}
\wt{f}(x_0) \leq \int{fd\mu_G} \leq \wt{f}(x_1),
\end{equation*}
since $\int{fd\mu_G}=\int{\wt{f}d\mu_G}$. \item
\emph{(Intermediate value theorem)} Finally we suppose (for a
contradiction) that $\chi_A \in A(\mathbb{T})$ so that by the
intermediate value theorem there is some $x \in \mathbb{T}$ such
that $\wt{\chi_A}(x)=\alpha$. Continuity ensures that there is an
open ball $x+B$ on which $\wt{\chi_A}$ is very close to $\alpha$,
and in particular, since $\alpha \in (0,1)$, on which
$\wt{\chi_A}$ only takes values in $(0,1)$. Since
$\|\chi_A-\wt{\chi_A}\|_{L^\infty(G)}=0$ and $\mu(x+B)>0$ it
follows that $\chi_A$ equals $\wt{\chi_A}$ for some point in
$x+B$, but this contradicts the fact that $\chi_A$ can only take
the values 0 or 1.
\end{enumerate}
If we try to transfer this argument to $G=\mathbb{Z}/p\mathbb{Z}$
it breaks down at the third stage when we apply the intermediate
value theorem. It is easy enough to remedy this and prove a
sensible discrete analogue of the intermediate value theorem; the
following, for example, is in \cite{BJGSVK}.
\begin{proposition}
\emph{(Discrete intermediate value theorem)}\label{divt} Suppose
that $p$ is prime number. Suppose that $f:\mathbb{Z}/p\mathbb{Z}
\rightarrow \mathbb{R}$ and that there is some non-zero $y \in
\mathbb{Z}/p\mathbb{Z}$ such that
\begin{equation*}
|f(x+y) -f(x)| \leq \epsilon\|f\|_\infty \textrm{ for all } x \in
\mathbb{Z}/p\mathbb{Z}.
\end{equation*}
Then there is some $x \in \mathbb{Z}/p\mathbb{Z}$ such that
\begin{equation*}
|f(x) - \int{fd\mu_{\mathbb{Z}/p\mathbb{Z}}}| \leq
2^{-1}\epsilon\|f\|_\infty.
\end{equation*}
\end{proposition}
Of course this has only moved the difficulty: to use this result
we need to replace the continuity in the first stage of our
argument with the sort of \emph{quantitative continuity} used in
this proposition.

It turns out that we already have a ready supply of functions which
are continuous in this new sense: Suppose that $f \in L^\infty(G)$
and $B(\Gamma,\delta)$ is a regular Bohr set. Then by Corollary
\ref{contlemcor2} we may pick $\delta' \gg \epsilon\delta/d$ such
that
\begin{equation*}
\|f\ast \beta-f\ast \beta(x)\|_{L^\infty(x+B(\Gamma,\delta'))} \leq
\epsilon\|f\|_{L^\infty(G)}.
\end{equation*}
Now if $\mu_G(B(\Gamma,\delta'))>p^{-1}$ then $B(\Gamma,\delta')$
has a non-identity element and hence the discrete intermediate
value theorem applies.

Essentially the same argument which shows that if $f \in A(G)$
then $\|f-\wt{f}\|_{L^\infty(G)}=0$ for some continuous function
$\wt{f}$, can be made quantitative to show that there is a regular
Bohr set $B(\Gamma,\delta)$ such that $\|f-f\ast
\beta\|_{L^\infty(G)}$ is small and, by our previous observations,
$f\ast \beta$ is quantitatively continuous.

To be concrete suppose that $G$ is a compact Abelian group, $f \in
A(G)$ and write $A_f:=\|f\|_{A(G)}\|f\|_{L^\infty(G)}^{-1}$. Then
there is a finite set of characters $\Gamma$ such that
\begin{equation*}
\sum_{\gamma \not\in \Gamma}{|\wh{f}(\gamma)|} \leq \epsilon
A_f^{-1} \|f\|_{A(G)}/3.
\end{equation*}
Pick $\delta \gg \epsilon A_f^{-1}$ such that
\begin{equation*}
B(\Gamma,\delta) \subset \{x \in G: |1-\gamma(x)| \leq \epsilon
A_f^{-1}/3 \textrm{ for all } \gamma \in \Gamma\},
\end{equation*}
and such that $\delta$ is regular for $\Gamma$ by Proposition
\ref{ubreg}. It is easy to see that
\begin{equation*}
|1-\wh{\beta}(\gamma)| \leq \epsilon A_f^{-1}/3 \textrm{ if } \gamma
\in \Gamma,
\end{equation*}
and it follows that
\begin{eqnarray*}
\|f-f \ast \beta\|_{L^\infty(G)} & \leq & \sum_{\gamma \in
\wh{G}}{|1-\wh{\beta}(\gamma)||\wh{f}(\gamma)|}\\ & \leq &
\sum_{\gamma \in \Gamma}{|1-\wh{\beta}(\gamma)||\wh{f}(\gamma)|} +
\sum_{\gamma \not \in \Gamma }{|1 - \wh{\beta}(\gamma)|
|\wh{f}(\gamma)|}\\ & \leq & (\epsilon A_f^{-1}/3)\sum_{\gamma \in
\wh{G}}{|\wh{f}(\gamma)|} + 2\sum_{\gamma \not \in
\Gamma}{|\wh{f}(\gamma)|}\\ & \leq & \epsilon A_f^{-1} \|f\|_{A(G)}
= \epsilon \|f\|_{L^\infty(G)}.
\end{eqnarray*}
In slightly formal language this has proved the following result.
\begin{theorem}\label{poor}
Suppose that $G$ is a compact Abelian group, $f \in A(G)$ and
$\epsilon \in (0,1]$. Write $A_f$ for the quantity
$\|f\|_{A(G)}\|f\|_{L^\infty(G)}^{-1}$. Then there is a Bohr set
$B(\Gamma,\delta)$ with
\begin{equation*}
d < \infty \textrm{ and } \delta^{-1} \ll \epsilon^{-1}A_f,
\end{equation*}
and a narrower Bohr set $B(\Gamma,\delta')$ with $\delta' \gg
\epsilon \delta/d$ such that
\begin{equation*}
\sup_{x \in G}{\|f\ast \beta-f \ast
\beta(x)\|_{L^\infty(x+B(\Gamma,\delta'))}} \leq \epsilon
\|f\|_{L^\infty(G)}
\end{equation*}
and
\begin{equation*} \sup_{x \in G}{\|f-f \ast
\beta\|_{L^\infty(x+B(\Gamma,\delta'))}} \leq \epsilon
\|f\|_{L^\infty(G)}.
\end{equation*}
\end{theorem}
Of course, as we observed before, this is only useful to us if
$B(\Gamma,\delta')$ contains a non-zero element. We can use Lemma
\ref{bohrsize} to estimate its size: If
$G=\mathbb{Z}/p\mathbb{Z}$, then $B(\Gamma,\delta')$ contains a
non-zero element if
\begin{equation*}
(c\epsilon^2A_f^{-1}/d)^d>p^{-1} \textrm{ for some absolute }c>0.
\end{equation*}
Unfortunately, because we have no control over $d$, we have no way
of ensuring this inequality. The content of this paper can be seen
as an effort to make this method work by getting control of $d$;
the main result is the following refinement of Theorem \ref{poor}.
\begin{theorem}\label{mainconttheorem}
Suppose that $G$ is a compact Abelian group, $f \in A(G)$ and
$\epsilon \in (0,1]$. Write $A_f$ for the quantity
$\|f\|_{A(G)}\|f\|_{\infty}^{-1}$. Then there is a Bohr set
$B(\Gamma,\delta)$ with
\begin{equation*}
d \ll \epsilon^{-2}A_f(1+ \log A_f)(1+ \log \epsilon^{-1}A_f)
\textrm{ and } \log \delta^{-1} \ll \epsilon^{-2}A_f (1+\log
\epsilon^{-1}A_f),
\end{equation*}
and a narrower Bohr set $B(\Gamma,\delta')$ with $\delta' \gg
\epsilon \delta/d$ such that
\begin{equation*}
\sup_{x \in G}{\|f\ast \beta-f \ast
\beta(x)\|_{L^\infty(x+B(\Gamma,\delta'))}} \leq \epsilon
\|f\|_{L^\infty(G)}
\end{equation*}
and
\begin{equation*} \sup_{x \in G}{\|f-f \ast
\beta\|_{L^2(x+B(\Gamma,\delta'))}} \leq \epsilon
\|f\|_{L^\infty(G)}.
\end{equation*}
\end{theorem}
Note that to gain control of $d$ we have had to sacrifice some
control of $\delta$ and of the error in approximating $f$ by $f
\ast \beta$.

There are now three remaining sections to the paper.
\begin{itemize}
\item \S\ref{modelargument} details our arguments in the model
setting of $\mathbb{F}_2^n$. \item \S\ref{mainarg} proves Theorem
\ref{mainconttheorem} following the outline of
\S\ref{modelargument} and using the tools of
\S\S\ref{localft}\verb.&.\ref{localstruct}. \item Finally
\S\ref{comp} completes the proof of Theorem \ref{maintheorem} and
concludes with some remarks and a conjecture.
\end{itemize}

\section{The argument in a model setting}\label{modelargument}

Many proofs of results for compact Abelian groups can be modelled
much more cleanly in $\mathbb{F}_2^n$, indeed we have already
partially seen this fact in \S\ref{auxiliarymeasure}, and so the
purpose of this section is to prove the following result, which is
the model version of the Theorem \ref{mainconttheorem}.
\begin{theorem}\label{f2ndetail}
Suppose that $G=\mathbb{F}_2^n$. Suppose that $f \in A(G)$ and
$\epsilon \in (0,1]$. Write $A_f$ for the quantity
$\|f\|_{A(G)}\|f\|_{L^\infty(G)}^{-1}$. Then there is a subspace $V$
of $G$ with
\begin{equation*}
\codim V \ll \epsilon^{-2}A_f (1+\log A_f)(1+ \log \epsilon^{-1}
A_f),
\end{equation*}
and
\begin{equation*}
\sup_{x \in G}{\|f-f\ast \mu_V\|_{L^2(x+V)}} \leq
\epsilon\|f\|_{L^\infty(G)}.
\end{equation*}
\end{theorem}
The first part of the conclusion of Theorem \ref{mainconttheorem}
is unnecessary since $f\ast \mu_V$ is constant on cosets of $V$
and hence
\begin{equation*}
\sup_{x \in G}{\|f\ast \mu_V-f \ast \mu_V(x)\|_{L^\infty(x+V)}}
=0.
\end{equation*}

\subsection{The basic quantitative argument}

We begin with an argument which proves the following weak version
of Theorem \ref{f2ndetail}; the argument will form the basis of
our proof of that theorem.
\begin{theorem}\label{f2nsimp}
Suppose that $G=\mathbb{F}_2^n$. Suppose that $f \in A(G)$ and
$\epsilon \in (0,1]$. Write $A_f$ for the quantity
$\|f\|_{A(G)}\|f\|_{L^\infty(G)}^{-1}$. Then there is a subspace $V$
of $G$ with
\begin{equation*}
\codim V \leq 2^3\epsilon^{-4}A_f^3,
\end{equation*}
and
\begin{equation*}
\sup_{x \in G}{\|f-f\ast \mu_V\|_{L^2(x+V)}} \leq
\epsilon\|f\|_{L^\infty(G)}.
\end{equation*}
\end{theorem}
The technique is iterative; the driving component is the following
lemma.
\begin{lemma}\label{itlem1}\emph{(Iteration lemma 1)}
Suppose that $G=\mathbb{F}_2^n$ and $\Gamma^\perp$ is an annihilator
in $G$. Suppose that $f \in A(G)$ and $\epsilon \in (0,1]$. Write
$A_f$ for the quantity $\|f\|_{A(G)}\|f\|_{L^\infty(G)}^{-1}$. Then
at least one of the following is true.
\begin{enumerate}
\item \emph{($f$ is close to a continuous
function)}\begin{equation*} \sup_{x \in G}{\|f-f \ast
\mu_{\Gamma^\perp}\|_{L^2(x+\Gamma^\perp)}} \leq \epsilon
\|f\|_{L^\infty(G)}.
\end{equation*} \item There is a set of characters $\Lambda$ with
$|\Lambda| \leq 2\epsilon^{-2}A_f^2$ such that
\begin{equation*}
\sum_{\gamma \in
(\Gamma\cup\Lambda)^{\perp\perp}}{|\wh{f}(\gamma)|} - \sum_{\gamma
\in \Gamma^{\perp\perp}}{|\wh{f}(\gamma)|} \geq
2^{-2}\epsilon^2\|f\|_{L^\infty(G)}.
\end{equation*}
\end{enumerate}
\end{lemma}
Essentially this says that if $f$ doesn't satisfy the conclusion
of Theorem \ref{f2nsimp} for some annihilator $\Gamma^\perp$ then
there is a smaller annihilator $\Gamma'^\perp$, which is not too
much smaller, that supports more $A(G)$-norm of $f$.

To control the size of $\Gamma'^\perp$ we use Proposition
\ref{localag} from \S\ref{localstruct}; in $\mathbb{F}_2^n$ its
statement is particularly simple:
\begin{proposition}\label{modlocag}\emph{(Model analogue of Proposition
\ref{localag})} Suppose that $G=\mathbb{F}_2^n$ and $\Gamma^\perp$
is an annihilator in $G$. Suppose that $f \in A(\Gamma^\perp)$ and
$\epsilon \in (0,1]$. Then there is a set $\Lambda$ of characters
with
\begin{equation*}
|\Lambda| \leq
\epsilon^{-1}\|f\|_{A(\Gamma^\perp)}\|f\|_{L^\infty(\Gamma^\perp)}^{-1}
\end{equation*}
such that
\begin{equation*}
\{\gamma \in \wh{G}:|\wh{fd\mu_{\Gamma^\perp}}(\gamma)| \geq
\epsilon \|f\|_{L^\infty(\Gamma^\perp)}\} \subset (\Gamma \cup
\Lambda)^{\perp\perp}.
\end{equation*}
\end{proposition}
\begin{proof}[Proof of Lemma \ref{itlem1}] Suppose that
\begin{equation*}
\sup_{x \in G}{\|f-f\ast \mu_{\Gamma^\perp}\|_{L^2(x+\Gamma^\perp)}}
> \epsilon \|f\|_{L^\infty(G)}.
\end{equation*}
Then there is some $x' \in G$ which, without loss of generality, is
equal to $0_G$ such that
\begin{equation}\label{leverage}
\|f-f\ast \mu_{\Gamma^\perp}\|_{L^2(x'+\Gamma^\perp)} \geq
\epsilon \|f\|_{L^\infty(G)}.
\end{equation}
For ease of notation write $g=f - f \ast \mu_{\Gamma^\perp}$, and
observe that $g$ satisfies the inequalities
\begin{equation*}
\|g\|_{A(G)} \leq \|f\|_{A(G)} \textrm{ and }
\|g\|_{L^{\infty}(\Gamma^\perp)} \leq 2\|f\|_{L^\infty(G)}.
\end{equation*}
To see the first of these note that
\begin{equation*}
\|g\|_{A(G)} = \sum_{\gamma \in
\wh{G}}{|1-\wh{\mu_{\Gamma^\perp}}(\gamma)||\wh{f}(\gamma)|} \leq
\sup_{\gamma \in \wh{G}}{|1-\wh{\mu_{\Gamma^\perp}}(\gamma)|}
\|f\|_{A(G)} \leq \|f\|_{A(G)},
\end{equation*}
and for the second
\begin{equation*}
\|g\|_{L^\infty(\Gamma^\perp)} \leq \|g\|_{L^\infty(G)} \leq
\|f\|_{L^\infty(G)} + \|f \ast \mu_{\Gamma^\perp}\|_{L^\infty(G)}
\leq 2 \|f\|_{L^\infty(G)}.
\end{equation*}

Returning to (\ref{leverage}) we have
\begin{eqnarray}
\nonumber \epsilon^2\|f\|_{L^\infty(G)}^2 & \leq &
\|g\|_{L^2(\Gamma^\perp)}^2\\ \nonumber & = & \sum_{\gamma \in
\wh{G}}{\wh{gd\mu_{\Gamma^\perp}}(\gamma)\overline{\wh{g}(\gamma)}}
\textrm{ by Plancherel's theorem,}\\ \label{hjy} & \leq
&\sum_{\gamma \in
\wh{G}}{|\wh{gd\mu_{\Gamma^\perp}}(\gamma)||\wh{g}(\gamma)|}
\textrm{ by the triangle inequality.}
\end{eqnarray}
The characters supporting large values of
$\wh{gd\mu_{\Gamma^\perp}}$ make the principal contribution to
this sum. Specifically put
\begin{equation*}
\mathcal{C}:=\{\gamma \in
\wh{G}:|\wh{gd\mu_{\Gamma^\perp}}(\gamma)| >
\epsilon'\|g\|_{L^\infty(\Gamma^\perp)}\},
\end{equation*}
where
\begin{equation*}
\epsilon':=2^{-1}\epsilon^2A_f^{-1}\|f\|_{L^\infty(G)}\|g\|_{L^\infty(\Gamma^\perp)}^{-1}.
\end{equation*}
Then
\begin{eqnarray*}
\sum_{\gamma \not \in
\mathcal{C}}{|\wh{gd\mu_{\Gamma^\perp}}(\gamma)||\wh{g}(\gamma)|}
& \leq & 2^{-1}\epsilon^2A_f^{-1}\|f\|_{L^\infty(G)}\sum_{\gamma
\not \in \mathcal{C}}{|\wh{g}(\gamma)|}\\ & \leq
&2^{-1}\epsilon^2A_f^{-1}\|f\|_{L^\infty(G)}\|g\|_{A(G)}\\ & \leq
& 2^{-1}\epsilon^2\|f\|_{L^\infty(G)}^2 \textrm{ since
$\|g\|_{A(G)} \leq \|f\|_{A(G)}$.}
\end{eqnarray*}
Substituting this into (\ref{hjy}) we conclude that
\begin{equation}\label{used}
\sum_{\gamma \in
\mathcal{C}}{|\wh{gd\mu_{\Gamma^\perp}}(\gamma)||\wh{g}(\gamma)|}
\geq 2^{-1} \epsilon^2 \|f\|_{L^\infty(G)}^2.
\end{equation}
Now certainly $|\wh{gd\mu_{\Gamma^\perp}}(\gamma)| \leq
2\|f\|_{L^\infty(G)}$ so that
\begin{equation*}
2^{-2} \epsilon^2 \|f\|_{L^\infty(G)} \leq \sum_{\gamma \in
\mathcal{C}}{|\wh{g}(\gamma)|}.
\end{equation*}
Since $\|g\|_{A(\Gamma^\perp)} \leq \|f\|_{A(G)}$ we may apply
Proposition \ref{modlocag} to $\mathcal{C}$ to get a set of
characters $\Lambda$ with
\begin{equation*}
|\Lambda| < (\epsilon')^{-1}
\|g\|_{A(\Gamma^\perp)}\|g\|_{L^\infty(\Gamma^\perp)}^{-1} \leq 2
\epsilon^{-2} A_f^2,
\end{equation*}
such that $\mathcal{C} \subset (\Gamma \cup
\Lambda)^{\perp\perp}$. The lemma follows.
\end{proof}
We are now in a position to iterate this and prove Theorem
\ref{f2nsimp}.
\begin{proof}
[Proof of Theorem \ref{f2nsimp}] We construct a sequence of
annihilators $\Gamma_k^\perp$ iteratively. Define
\begin{equation*}
L_k:=\sum_{\gamma \in \Gamma_k^{\perp\perp}}{|\wh{f}(\gamma)|},
\end{equation*}
and initiate the iteration with $\Gamma_0:=\{0_{\wh{G}}\}$.

Suppose that we are at stage $k$ of the iteration. Apply the
iteration lemma (Lemma \ref{itlem1}). If we are in the first case
of the lemma then put $V=\Gamma_k^\perp$ and terminate; if not
then we get a set of characters $\Lambda$ and put
$\Gamma_{k+1}=\Gamma_k\cup\Lambda$. It follows from the properties
of $\Lambda$ that
\begin{equation*}
|\Gamma_{k+1}| \leq |\Gamma_k| + 2\epsilon^{-2}A_f^2 \textrm{ and
} L_{k+1}-L_k \geq 2^{-2}\epsilon^2\|f\|_{L^\infty(G)}.
\end{equation*}
By induction we have that after $k$ iterations
\begin{equation*}
|\Gamma_k| \leq k.2\epsilon^{-2}A_f^2 \textrm{ and } L_k \geq
k.2^{-2}\epsilon^2\|f\|_{L^\infty(G)}.
\end{equation*}
Since $L_k \leq \|f\|_{A(G)}$ we conclude that the iteration
terminates and
\begin{equation*}
|\Gamma_k| \leq 2^3\epsilon^{-4}A_f^3.
\end{equation*}
The theorem follows.
\end{proof}

\subsection{Refining the basic argument: the proof of Theorem \ref{f2ndetail}}

To achieve the result in Theorem \ref{f2ndetail} we make two
important improvements to the iteration lemma (Lemma \ref{itlem1})
of the previous argument.
\begin{itemize}
\item \emph{(Dyadic decomposition)} Our first improvement is the
observation that having derived (\ref{used}):
\begin{equation*}
\sum_{\gamma \in
\mathcal{C}}{|\wh{gd\mu_{\Gamma^\perp}}(\gamma)||\wh{g}(\gamma)|}
\geq 2^{-1}\epsilon^2\|f\|_{L^\infty(G)}^2,
\end{equation*}
we can do something better that simply adding all the characters in
$\mathcal{C}$ to $\Gamma$. Partition the characters in $\mathcal{C}$
by dyadically decomposing the range of values of
$|\wh{gd\mu_{\Gamma^\perp}}|$ and pick the characters in a dyadic
class contributing maximal mass to (\ref{used}). The $A(G)$-norm of
$f$ supported on this class is more closely related to the size of
$\mathcal{C}$ which yields the first improvement. \item
\emph{(Structure theorem for the Fourier spectrum)} The second
improvement replaces the application of Proposition \ref{modlocag}
with the stronger Proposition \ref{f2agchang}, which in the model
setting has the following simpler statement.
\begin{proposition}\label{modlocag2}\emph{(Model analogue of
Proposition \ref{f2agchang})} Suppose that $G=\mathbb{F}_2^n$ and
$\Gamma^\perp$ is an annihilator in $G$. Suppose that $f \in
A(\Gamma^\perp)$ and $\epsilon \in (0,1]$. Write $A_f$ for the
quantity
$\|f\|_{A(\Gamma^\perp)}\|f\|_{L^\infty(\Gamma^\perp)}^{-1}$. Then
there is a set $\Lambda$ of characters with $|\Lambda| \ll
\epsilon^{-1}(1+\log A_f)$ such that
\begin{equation*}
\{\gamma \in \wh{G}:|\wh{fd\mu_{\Gamma^\perp}}(\gamma)| \geq
\epsilon \|f\|_{L^\infty(\Gamma^\perp)}\} \subset (\Gamma \cup
\Lambda)^{\perp\perp}.
\end{equation*}
\end{proposition}
\end{itemize}
By implementing these two refinements we prove the following
iteration lemma.
\begin{lemma}\emph{(Iteration lemma 2)}\label{itlem2}
Suppose that $G=\mathbb{F}_2^n$ and $\Gamma^\perp$ is an annihilator
in $G$. Suppose that $f \in A(G)$ and $\epsilon \in (0,1]$. Write
$A_f$ for the quantity $\|f\|_{A(G)}\|f\|_{L^\infty(G)}^{-1}$. Then
at least one of the following is true.
\begin{enumerate}
\item \emph{($f$ is close to a continuous
function)}\begin{equation*} \sup_{x \in G}{\|f-f \ast
\mu_{\Gamma^\perp}\|_{L^2(x+\Gamma^\perp)}} \leq \epsilon
\|f\|_{L^\infty(G)}.
\end{equation*} \item There is a set of characters $\Lambda$ and a non-negative integer $s$ with
$|\Lambda| \ll 2^s(1+\log A_f)$ such that
\begin{equation*}
\sum_{\gamma \in (\Gamma\cup
\Lambda)^{\perp\perp}}{|\wh{f}(\gamma)|} - \sum_{\gamma \in
\Gamma^{\perp\perp}}{|\wh{f}(\gamma)|} \gg
\frac{2^s\epsilon^2\|f\|_{L^\infty(G)}}{1+\log \epsilon^{-1}A_f}.
\end{equation*}
\end{enumerate}
\end{lemma}
\begin{proof}
We proceed as in the proof of Lemma \ref{itlem1} up to the point
where we conclude that
\begin{equation*}
\sum_{\gamma \in
\mathcal{C}}{|\wh{gd\mu_{\Gamma^\perp}}(\gamma)||\wh{g}(\gamma)|}
\geq 2^{-1}\epsilon^2 \|f\|_{L^\infty(G)}^2 \textrm{
(\ref{used}).}
\end{equation*}
Write $I_s:=(2^{-s}\|f\|_{L^\infty(G)},
2^{-(s-1)}\|f\|_{L^\infty(G)}]$ and partition $\mathcal{C}$ into
the sets
\begin{equation*}
\mathcal{C}_s:=\{\gamma \in
\mathcal{C}:|\wh{gd\mu_{\Gamma^\perp}}(\gamma)| \in I_s\} \textrm{
for } 0 \leq s \leq 2+ \log_2 \epsilon^{-2}A_f.
\end{equation*}
Note that $\{\mathcal{C}_s: 0 \leq s \leq 2 +
\log_2\epsilon^{-2}A_f\}$ covers $\mathcal{C}$ since
\begin{equation*}
\sup_{\gamma \in \mathcal{C}}{|\wh{gd\mu_{\Gamma^\perp}}(\gamma)|}
\leq \sup_{\gamma \in \wh{G}}{|\wh{gd\mu_{\Gamma^\perp}}(\gamma)|}
\leq \|g\|_{L^\infty(\Gamma^\perp)} \leq 2\|f\|_{L^\infty(G)}
\end{equation*}
and
\begin{equation*}
\inf_{\gamma \in\mathcal{C}}{|\wh{gd\mu_{\Gamma^\perp}}(\gamma)|}
> 2^{-1}\epsilon^2A_f^{-1}\|f\|_{L^\infty(G)},
\end{equation*}
so that (\ref{used}) may be rewritten to yield
\begin{equation*}
\sum_{s=0}^{2+ \log_2 \epsilon^{-2}A_f}{\sum_{\gamma \in
\mathcal{C}_s}{|\wh{gd\mu_{\Gamma^\perp}}(\gamma)||\wh{g}(\gamma)|}}
\geq 2^{-1} \epsilon^2 \|f\|_{L^\infty(G)}^2.
\end{equation*}
It follows by the pigeonhole principle that there is some $s$ for
which
\begin{equation*}
\sum_{\gamma \in
\mathcal{C}_s}{|\wh{gd\mu_{\Gamma^\perp}}(\gamma)||\wh{g}(\gamma)|}
\gg \frac{\epsilon^2\|f\|_{L^\infty(G)}^2}{ 1+\log
\epsilon^{-1}A_f},
\end{equation*}
and since $|\wh{gd\mu_{\Gamma^\perp}}(\gamma)| \leq
2^{-(s-1)}\|f\|_{L^\infty(G)}$ if $\gamma \in \mathcal{C}_s$ we
get
\begin{equation*}
\sum_{\gamma \in \mathcal{C}_s}{|\wh{g}(\gamma)|} \gg
\frac{2^s\epsilon^2\|f\|_{L^\infty(G)}}{1+\log \epsilon^{-1}A_f}.
\end{equation*}

Now
\begin{equation*}
\mathcal{C}_s \subset \{ \gamma : |\wh{gd\beta'}(\gamma)| >
(2^{-s}\|f\|_{L^\infty(G)}\|g\|_{L^\infty(\Gamma^\perp)}^{-1})\|g\|_{L^\infty(\Gamma^\perp)}\},
\end{equation*}
and since $\|g\|_{A(\Gamma^\perp)} \leq \|f\|_{A(G)}$ we may apply
Proposition \ref{modlocag2} to get a set of characters $\Lambda$
such that $\mathcal{C}_s \subset (\Gamma \cup
\Lambda)^{\perp\perp}$. Moreover $|\Lambda|$ satisfies
\begin{eqnarray*}
|\Lambda| & \ll &
2^s\|f\|_{L^\infty(G)}^{-1}\|g\|_{L^\infty(\Gamma^\perp)}(1+ \log
\|g\|_{A(\Gamma^\perp)}\|g\|_{L^\infty(\Gamma^\perp)}^{-1})\\
& \ll & 2^s
\|f\|_{L^\infty(G)}^{-1}\|g\|_{L^\infty(\Gamma^\perp)}(1+ \log
\|f\|_{A(G)}\|g\|_{L^\infty(\Gamma^\perp)}^{-1})\textrm{ since
}\|g\|_{A(\Gamma^\perp)} \leq \|f\|_{A(G)}\\ & \ll & 2^s
\|f\|_{L^\infty(G)}^{-1}\|g\|_{L^\infty(\Gamma^\perp)}(1+ \log
A_f\|f\|_{L^\infty(G)}\|g\|_{L^\infty(\Gamma^\perp)}^{-1}).
\end{eqnarray*}
So, writing $X$ for
$\|f\|_{L^\infty(G)}^{-1}\|g\|_{L^\infty(\Gamma^\perp)}$ we have
\begin{equation*}
|\Lambda| \ll 2^sX(1+\log A_f X^{-1}),
\end{equation*}
but $\|g\|_{L^\infty(\Gamma^\perp)} \leq 2\|f\|_{L^\infty(G)}$ so $X
\leq 2$ and therefore
\begin{equation*}
|\Lambda| \ll 2^s \sup_{X' \in (0,2]}{X'(1+\log A_f X'^{-1})}\ll
2^s(1+\log A_f).
\end{equation*}
The lemma follows.
\end{proof}
Iterating this in the same way as before yields Theorem
\ref{f2ndetail}.

\section{The proof of Theorem \ref{mainconttheorem}}\label{mainarg}

We begin by extending the second iteration lemma (Lemma
\ref{itlem2}) from the model setting to that of the general compact
Abelian group.
\begin{lemma}\label{itlemcont}
Suppose that $G$ is a compact Abelian group and $B(\Gamma,\delta)$ a
regular Bohr set. Suppose that $f \in A(G)$ and $\epsilon \in
(0,1]$. Write $A_f$ for the quantity
$\|f\|_{A(G)}\|f\|_{L^\infty(G)}^{-1}$. Then at least one of the
following is true.
\begin{enumerate}
\item \emph{($f$ is close to a continuous function)} There is a
Bohr set $B(\Gamma,\delta')$ with $\delta' \gg \epsilon \delta/d$
such that
\begin{equation*}
\sup_{x \in G}{\|f\ast \beta-f \ast
\beta(x)\|_{L^\infty(x+B(\Gamma,\delta'))}} \leq \epsilon
\|f\|_{L^\infty(G)}
\end{equation*}
and
\begin{equation*} \sup_{x \in G}{\|f-f \ast
\beta\|_{L^2(x+B(\Gamma,\delta'))}} \leq \epsilon
\|f\|_{L^\infty(G)}.
\end{equation*}
\item For all $\eta \in (0,1]$ there is a set of characters
$\Lambda$, a $\delta''\in (0,1]$ and a non-negative integer $s$
with
\begin{equation*}
|\Lambda| \ll 2^s(1+\log A_f) \textrm{ and } \delta'' \gg
\epsilon^5A_f^{-4}\eta\delta /d^{3},
\end{equation*}
such that
\begin{equation*}
\sum_{\gamma \in
\mathcal{L}}{|1-\wh{\beta}(\gamma)||\wh{f}(\gamma)|} \gg
\frac{2^s\epsilon^2\|f\|_{L^\infty(G)}}{\min\{2^s,1+\log
\epsilon^{-1}A_f\}}
\end{equation*}
where
\begin{equation*}
\mathcal{L}:=\{\gamma:|1-\gamma(x)| \leq \eta \textrm{ for all } x
\in B(\Gamma\cup\Lambda,\delta'')\}.
\end{equation*}
\end{enumerate}
\end{lemma}
\begin{proof}
Choosing $\delta'$ is easy: By Corollary \ref{contlemcor2} and
Proposition \ref{ubreg} there is a $\delta' \gg \delta\epsilon /d$
regular for $\Gamma$ such that
\begin{equation*}
\sup_{x \in G}{\|f\ast \beta-f \ast
\beta(x)\|_{L^\infty(x+B(\Gamma,\delta'))}} \leq \epsilon
\|f\|_{L^\infty(G)}.
\end{equation*}
Now, suppose that
\begin{equation*}
\sup_{x \in G}{\|f-f\ast \beta\|_{L^2(x+B(\Gamma,\delta'))}} >
\epsilon \|f\|_{L^\infty(G)}.
\end{equation*}
It follows that there is some $x' \in G$ which, without loss of
generality, is equal to $0_G$ such that
\begin{equation}\label{leverageprop}
\|f-f\ast \beta\|_{L^2(x'+B(\Gamma,\delta'))} \geq \epsilon
\|f\|_{L^\infty(G)}.
\end{equation}
For ease of notation write $g=f - f \ast \beta$, and observe that
$g$ satisfies the inequalities
\begin{equation*}
\|g\|_{A(G)} \leq 2\|f\|_{A(G)} \textrm{ and }
\|g\|_{L^{\infty}(B(\Gamma,\delta'))} \leq 2\|f\|_{L^\infty(G)}.
\end{equation*}
To see the first of these note that
\begin{equation*}
\|g\|_{A(G)} = \sum_{\gamma \in
\wh{G}}{|1-\wh{\beta}(\gamma)||\wh{f}(\gamma)|} \leq \sup_{\gamma
\in \wh{G}}{|1-\wh{\beta}(\gamma)|} \|f\|_{A(G)} \leq 2
\|f\|_{A(G)},
\end{equation*}
and for the second
\begin{equation*}
\|g\|_{L^\infty(B(\Gamma,\delta'))} \leq \|g\|_{L^\infty(G)} \leq
\|f\|_{L^\infty(G)} + \|f \ast \beta\|_{L^\infty(G)} \leq 2
\|f\|_{L^\infty(G)}.
\end{equation*}

Returning to (\ref{leverageprop}) we may apply Plancherel's
theorem and then the triangle inequality to give us a Fourier
statement:
\begin{equation}\label{hjyprop}
\sum_{\gamma \in \wh{G}}{|\wh{gd\beta'}(\gamma)||\wh{g}(\gamma)|}
\geq \epsilon^2 \|f\|_{L^\infty(G)}^2.
\end{equation}
The characters supporting large values of $\wh{gd\beta'}$ make the
principal contribution to this sum. Specifically put
\begin{equation*}
\mathcal{C}:=\{\gamma \in \wh{G}:|\wh{gd\beta'}(\gamma)| >
\epsilon'\|g\|_{L^\infty(B(\Gamma,\delta'))}\},
\end{equation*}
where
\begin{equation*}
\epsilon':=2^{-2}\epsilon^2A_f^{-1}\|f\|_{L^\infty(G)}\|g\|_{L^\infty(B(\Gamma,\delta'))}^{-1}.
\end{equation*}
Then
\begin{eqnarray*}
\sum_{\gamma \not \in
\mathcal{C}}{|\wh{gd\beta'}(\gamma)||\wh{g}(\gamma)|} & \leq &
2^{-2}\epsilon^2A_f^{-1}\|f\|_{L^\infty(G)}\sum_{\gamma \not \in
\mathcal{C}}{|\wh{g}(\gamma)|}\\ & \leq
&2^{-2}\epsilon^2A_f^{-1}\|f\|_{L^\infty(G)}\|g\|_{A(G)}\\ & \leq
& 2^{-1}\epsilon^2\|f\|_{L^\infty(G)}^2 \textrm{ since
$\|g\|_{A(G)} \leq 2\|f\|_{A(G)}$.}
\end{eqnarray*}
Substituting this into (\ref{hjyprop}) we conclude that
\begin{equation}\label{usedprop}
\sum_{\gamma \in
\mathcal{C}}{|\wh{gd\beta'}(\gamma)||\wh{g}(\gamma)|} \geq 2^{-1}
\epsilon^2 \|f\|_{L^\infty(G)}^2 .
\end{equation}
Write $I_s:=(2^{-s}\|f\|_{L^\infty(G)},
2^{-(s-1)}\|f\|_{L^\infty(G)}]$ and partition $\mathcal{C}$ into
the sets
\begin{equation*}
\mathcal{C}_s:=\{\gamma \in \mathcal{C}:|\wh{gd\beta'}(\gamma)| \in
I_s\} \textrm{ for } 0 \leq s \leq 3+ \log_2 \epsilon^{-2}A_f.
\end{equation*}
Notice that $\{\mathcal{C}_s: 0 \leq s \leq 3 +
\log_2\epsilon^{-2}A_f\}$ covers $\mathcal{C}$ since
\begin{equation*}
\sup_{\gamma \in \mathcal{C}}{|\wh{gd\beta'}(\gamma)|} \leq
\sup_{\gamma \in \wh{G}}{|\wh{gd\beta'}(\gamma)|} \leq
\|g\|_{L^\infty(B(\Gamma,\delta'))} \leq 2\|f\|_{L^\infty(G)}
\end{equation*}
and
\begin{equation*}
\inf_{\gamma \in\mathcal{C}}{|\wh{gd\beta'}(\gamma)|} >
2^{-2}\epsilon^2A_f^{-1}\|f\|_{L^\infty(G)},
\end{equation*}
so that (\ref{usedprop}) may be rewritten to yield
\begin{equation*}
\sum_{s=0}^{3+ \log_2 \epsilon^{-2}A_f}{\sum_{\gamma \in
\mathcal{C}_s}{|\wh{gd\beta'}(\gamma)||\wh{g}(\gamma)|}} \geq 2^{-1}
\epsilon^2 \|f\|_{L^\infty(G)}^2.
\end{equation*}
Writing $S':=\{s \in \mathbb{N}_0: 2^s \leq 3+ \log_2
\epsilon^{-2}A_f\}$ and $S'':=\{s \in \mathbb{N}_0: 2^s > 3+ \log_2
\epsilon^{-2}A_f\}$ it follows that either
\begin{equation*}
\sum_{s \in S'}{2^{-s}.2^s\sum_{\gamma \in
\mathcal{C}_s}{|\wh{gd\beta'}(\gamma)||\wh{g}(\gamma)|}} \geq 2^{-2}
\epsilon^2 \|f\|_{L^\infty(G)}^2
\end{equation*}
or
\begin{equation*}
\sum_{s \in S''}{\sum_{\gamma \in
\mathcal{C}_s}{|\wh{gd\beta'}(\gamma)||\wh{g}(\gamma)|}} \geq 2^{-2}
\epsilon^2 \|f\|_{L^\infty(G)}^2.
\end{equation*}
By the pigeonhole principle that there is some $s$ for which
\begin{equation*}
\sum_{\gamma \in
\mathcal{C}_s}{|\wh{gd\beta'}(\gamma)||\wh{g}(\gamma)|} \gg
\frac{\epsilon^2\|f\|_{L^\infty(G)}^2}{ 1+\log \epsilon^{-1}A_f}
\textrm{ if } 2^s > 3+\log_2 \epsilon^{-2}A_f
\end{equation*}
and
\begin{equation*}
\sum_{\gamma \in
\mathcal{C}_s}{|\wh{gd\beta'}(\gamma)||\wh{g}(\gamma)|} \gg
\frac{\epsilon^2\|f\|_{L^\infty(G)}^2}{ 2^s} \textrm{ if } 2^s \leq
3+\log_2 \epsilon^{-2}A_f.
\end{equation*}
i.e. there is some $s$ such that
\begin{equation*}
\sum_{\gamma \in
\mathcal{C}_s}{|\wh{gd\beta'}(\gamma)||\wh{g}(\gamma)|} \gg
\frac{\epsilon^2\|f\|_{L^\infty(G)}^2}{
\min\{2^s,1+\log\epsilon^{-1}A_f\}}.
\end{equation*}
Since $|\wh{gd\beta'}(\gamma)| \leq 2^{-(s-1)}\|f\|_{L^\infty(G)}$
if $\gamma \in \mathcal{C}_s$ we get
\begin{equation*}
\sum_{\gamma \in \mathcal{C}_s}{|\wh{g}(\gamma)|} \gg
\frac{2^s\epsilon^2\|f\|_{L^\infty(G)}}{\min\{2^s,1+\log
\epsilon^{-1}A_f\}}.
\end{equation*}
Now
\begin{equation*}
\mathcal{C}_s \subset \{ \gamma : |\wh{gd\beta'}(\gamma)| \geq
(2^{-s}\|f\|_{L^\infty(G)}\|g\|_{L^\infty(B(\Gamma,\delta'))}^{-1})\|g\|_{L^\infty(B(\Gamma,\delta'))}\},
\end{equation*}
and $g \in A(G)$ so we may apply Proposition \ref{localagchang} to
get a set of characters $\Lambda$ and a $\delta''$ regular for
$\Gamma \cup \Lambda$ such that
\begin{equation*}
\mathcal{C}_s \subset \{ \gamma: |1-\gamma(x)| \leq \eta \textrm{
for all }x \in B(\Gamma\cup \Lambda,\delta'')\}.
\end{equation*}
Moreover $|\Lambda|$ satisfies
\begin{eqnarray*} |\Lambda| & \ll &
2^s\|f\|_{L^\infty(G)}^{-1}\|g\|_{L^\infty(B(\Gamma,\delta'))}(1+
\log
\|g\|_{A(G)}\|g\|_{L^\infty(B(\Gamma,\delta'))}^{-1})\\
& \ll & 2^s
\|f\|_{L^\infty(G)}^{-1}\|g\|_{L^\infty(B(\Gamma,\delta'))}(1+ \log
2\|f\|_{A(G)}\|g\|_{L^\infty(B(\Gamma,\delta'))}^{-1})
\end{eqnarray*}
since $\|g\|_{A(G)} \leq 2\|f\|_{A(G)}$, so
\begin{equation*}
|\Lambda| \ll 2^s
\|f\|_{L^\infty(G)}^{-1}\|g\|_{L^\infty(B(\Gamma,\delta'))}(1+ \log
2A_f\|f\|_{L^\infty(G)}\|g\|_{L^\infty(B(\Gamma,\delta'))}^{-1}).
\end{equation*}
So, writing $X$ for
$\|f\|_{L^\infty(G)}^{-1}\|g\|_{L^\infty(B(\Gamma,\delta'))}$ we
have
\begin{equation*}
|\Lambda| \ll 2^sX(1+\log 2A_f X^{-1}),
\end{equation*}
but $\|g\|_{L^\infty(B(\Gamma,\delta'))} \leq 2\|f\|_{L^\infty(G)}$
so $X \leq 2$ and therefore
\begin{equation*}
|\Lambda| \ll 2^s \sup_{X' \in (0,2]}{X'(1+\log 2A_f X'^{-1})}\ll
2^s(1+\log A_f).
\end{equation*}
Furthermore $\delta''$ satisfies
\begin{eqnarray*}
\delta'' &\gg &
2^{-2s}\|f\|_{L^\infty(G)}^2\|g\|_{L^\infty(B(\Gamma,\delta'))}^{-2}\eta
\delta'/d^2(1+\log \|g\|_{A(G)}\|g\|_{L^\infty(B(\Gamma,\delta'))}^{-1})\\
& \gg & 2^{-2s}\|f\|_{L^\infty(G)}^2\eta
\delta'/d^2\|g\|_{A(G)}^2\\ & \gg &
\epsilon^4A_f^{-2}\|f\|_{L^\infty(G)}^2\eta
\delta'/d^2\|g\|_{A(G)}^2 \textrm{ since }2^{2s}
\leq 2^4\epsilon^{-4}A_f^2\\
 & \gg & \epsilon^4A_f^{-4}\eta\delta'/d^2 \textrm{ since
 }\|g\|_{A(G)} \leq 2\|f\|_{A(G)}.
\end{eqnarray*}
The lemma follows.
\end{proof}
We are now in a position to iterate this lemma.
\begin{proof}
[Proof of Theorem \ref{mainconttheorem}] Fix $\eta$ to be
optimized at the end of the argument. We construct a sequence of
regular Bohr sets $B(\Gamma_k,\delta_k)$ iteratively using Lemma
\ref{itlemcont}. Put
\begin{equation*}
\mathcal{L}_k:=\{\gamma:|1-\gamma(x)| \leq \eta \textrm{ for all
}x \in B(\Gamma_k,\delta_k)\}
\end{equation*}
and
\begin{equation*}
d_k:=|\Gamma_k| \textrm{ and } L_k:=\sum_{\gamma \in
\mathcal{L}_k}{|\wh{f}(\gamma)|}.
\end{equation*}
We initialize the iteration with $\Gamma_0:=\{0_{\wh{G}}\}$ and
$\delta_0 \gg 1$ regular for $\Gamma_0$, chosen so by Proposition
\ref{ubreg}.

Suppose that we are at stage $k$. Apply the iteration lemma (Lemma
\ref{itlemcont}) to $f$ and the regular Bohr set
$B(\Gamma_k,\delta_k)$. If we are in the first case terminate with
the desired conclusion; if not then we get a set of characters
$\Lambda$, a $\delta'' \in (0,1]$ and an integer $s$. Let
$\Gamma_{k+1}=\Gamma_k\cup \Lambda$, pick $\delta_{k+1} \in
(\delta''/2,\delta'']$ regular for $\Gamma_{k+1}$ by Proposition
\ref{ubreg}, and let $s_{k+1}=s$. We are given that
\begin{equation*}
d_{k+1} - d_k \ll 2^{s_{k+1}}(1+\log A_f) \textrm{ and }
\delta_{k+1} \gg \epsilon^5A_f^{-4} \eta \delta_k/d_k^3,
\end{equation*}
and furthermore
\begin{equation*}
2(L_{k+1} -L_k) + \eta L_k \gg
\frac{2^{s_{k+1}}\epsilon^2\|f\|_{L^\infty(G)}}{\min\{2^{s_{k+1}},1+\log
\epsilon^{-1}A_f\}}.
\end{equation*}
Since $L_k \leq \|f\|_{A(G)}$ and $s_k \geq 0$ it follows that we
can pick $\eta \gg \epsilon^3A_f^{-2}$ (independently of $k$) such
that
\begin{equation*}
L_{k+1} - L_k \gg
\frac{2^{s_{k+1}}\epsilon^2\|f\|_{L^\infty(G)}}{\min\{2^{s_{k+1}},1+\log
\epsilon^{-1}A_f\}}.
\end{equation*}
Hence by induction we have
\begin{equation*}
L_k \gg
\epsilon^2\|f\|_{L^\infty(G)}\sum_{l=1}^k{\frac{2^{s_{l}}}{\min\{2^{s_l},1+\log
\epsilon^{-1}A_f\}}} \textrm{ and } d_k \ll
\sum_{l=1}^k{2^{s_{l}}(1+\log A_f )}.
\end{equation*}
Again since $s_k \geq 0$ it follows that the iteration terminates.
Hence we have
\begin{equation*}
\sum_{l=1}^k{2^{s_l}} \ll L_k\epsilon^{-2} \|f\|_{L^\infty(G)}^{-1}
(1+\log \epsilon^{-1} A_f) \ll \epsilon^{-2} A_f (1+\log
\epsilon^{-1} A_f),
\end{equation*}
since $\|f\|_{A(G)} \geq L_k$. It follows that
\begin{equation*}
d_k \ll \epsilon^{-2}A_f (1+\log A_f)(1+ \log \epsilon^{-1}A_f).
\end{equation*}
The bound on $\eta$ and $d_k$ gives us
\begin{equation*}
\delta_{k+1} \gg \epsilon^{17}A_f^{-15}\delta_k,
\end{equation*}
and hence
\begin{eqnarray*}
\log \delta_k^{-1} & \ll & k (1+\log \epsilon^{-1} A_f)\\ & \ll &
\sum_{l=1}^k{\frac{2^{s_k}}{\min\{2^{s_k},1+\log
\epsilon^{-1}A_f\}}}(1+\log \epsilon^{-1} A_f)\\ & \ll &
\epsilon^{-2}A_f (1+\log \epsilon^{-1} A_f).
\end{eqnarray*}
The result follows.
\end{proof}

\section{The proof of Theorem \ref{maintheorem} and concluding
remarks} \label{comp}

Having proved Theorem \ref{mainconttheorem} it is essentially a
formality to carry out the rest of the argument detailed in
\S\ref{tfr}.

\begin{proof}
[Proof of Theorem \ref{maintheorem}] Write $G$ for
$\mathbb{Z}/p\mathbb{Z}$ and $\alpha:=\mu_G(A)=|A|/p$. We apply
Theorem \ref{mainconttheorem} to $f=\chi_A$ with $\epsilon
=2^{-2}\alpha(1-\alpha)$. This gives a Bohr set $B(\Gamma,\delta)$
with
\begin{equation*}
d \ll_{\alpha} \|\chi_A\|_{A(G)} (1+\log \|\chi_A\|_{A(G)})^2
\end{equation*}
and
\begin{equation*}
\log \delta^{-1} \ll_{\alpha} \|\chi_A\|_{A(G)} (1+\log
\|\chi_A\|_{A(G)}),
\end{equation*}
 and a narrower Bohr set $B(\Gamma,\delta')$ with $\delta'
\gg_\alpha \delta/d$ such that
\begin{equation*}
\sup_{x \in G}{\|\chi_A\ast \beta-\chi_A \ast
\beta(x)\|_{L^\infty(x+B(\Gamma,\delta'))}} \leq
2^{-2}\alpha(1-\alpha)
\end{equation*}
and
\begin{equation}\label{tge} \sup_{x \in G}{\|\chi_A-\chi_A \ast
\beta\|_{L^2(x+B(\Gamma,\delta'))}} \leq 2^{-2}\alpha(1-\alpha).
\end{equation}

Suppose that $\mu_G(B(\Gamma,\delta'))>p^{-1}$. Then there is a
non-zero $y \in B(\Gamma,\delta')$, and such a $y$ has the
property that $|\chi_A\ast \beta(x+y) - \chi_A \ast \beta(x)| \leq
2^{-2}\alpha(1-\alpha)$ for all $x \in G$. It follows that we may
apply the discrete intermediate value theorem (Proposition
\ref{divt}) to $\chi_A \ast \beta$ and conclude that there is some
$x \in G$ such that
\begin{equation*}
|\chi_A \ast \beta(x) - \alpha| \leq 2^{-3}\alpha(1-\alpha).
\end{equation*}
Furthermore (\ref{tge}) ensures that there is some $x' \in
x+B(\Gamma,\delta')$ such that
\begin{equation*}
|\chi_A(x') - \chi_A \ast \beta(x')| \leq 2^{-2}\alpha(1-\alpha),
\end{equation*}
and this gives
\begin{eqnarray*}
|\chi_A(x')-\alpha| & \leq & |\chi_A(x') - \chi_A \ast \beta(x')|\\
& & + |\chi_A\ast \beta(x') - \chi_A \ast \beta(x)|\\ & & + |\chi_A
\ast \beta(x) - \alpha|\\ &  \leq &\alpha(1-\alpha).
\end{eqnarray*}
This contradicts the fact that $\chi_A(x') \in \{0,1\}$, and hence
$\mu_G(B(\Gamma,\delta')) \leq p^{-1}$. Lemma \ref{bohrsize} then
lets us infer that $d (1+\log \delta'^{-1}) \gg \log p$ from which,
on inserting the bounds on $d$ and $\delta'$, the result follows.
\end{proof}

In \cite{BJGSVK} Green and Konyagin essentially prove a version of
Theorem \ref{mainconttheorem} with different bounds.
\begin{theorem}\label{green-kony}
Suppose that $G$ is a compact Abelian group, $f \in A(G)$ and
$\epsilon \in (0,1]$. Write $A_f$ for the quantity
$\|f\|_{A(G)}\|f\|_{\infty}^{-1}$. Then there is a Bohr set
$B(\Gamma,\delta)$ with
\begin{equation*}
|\Gamma| \ll \epsilon^{-2}A_f^2 \textrm{ and } \log \delta^{-1} \ll
\epsilon^{-1}A_f (1+\log \epsilon^{-1}A_f),
\end{equation*}
and a narrower Bohr set $B(\Gamma,\delta')$ with $\delta' \gg
\epsilon \delta/d$ such that
\begin{equation*}
\sup_{x \in G}{\|f\ast \beta-f \ast
\beta(x)\|_{L^\infty(x+B(\Gamma,\delta'))}} \leq \epsilon
\|f\|_{L^\infty(G)}
\end{equation*}
and
\begin{equation*} \sup_{x \in G}{\|f-f \ast
\beta\|_{L^2(x+B(\Gamma,\delta'))}} \leq \epsilon
\|f\|_{L^\infty(G)}.
\end{equation*}
\end{theorem}
The crucial difference between our proof of Theorem
\ref{mainconttheorem} and their proof of Theorem \ref{green-kony}
is that in their iteration lemma they find only a few characters
at which $\wh{f}$ is large, whereas we find all characters at
which $\wh{f}$ is large. Their approach leads to superior bounds
in the basic version of their argument, however it prevents them
from using a tool such as Proposition \ref{localagchang}, which is
where our argument gains its edge.

In both the arguments of Green and Konyagin and of this paper the
width of the Bohr set which one eventually finds narrows
exponentially with the number of times one has to use the
(appropriate) iteration lemma. Green and Konyagin employ a neat
trick to reduce this - the natural version of their argument has
$\log \delta^{-1} \ll \epsilon^{-2} A_f^2(1+\log \epsilon^{-1}A_f)$
- which leads to the superior $\epsilon$-dependence for $\log
\delta^{-1}$ in Theorem \ref{green-kony}. It is possible to add
their trick to our argument and hence improve the
$\epsilon$-dependence of $\log \delta^{-1}$ in Theorem
\ref{mainconttheorem} too, however this would have no effect on our
application.

The model setting of $\mathbb{F}_2^n$, has been used extensively
in this paper to make our results clearer; the paper \cite{BJGFFM}
of Green serves as a good survey of other problems where it has
found similar uses. While the primary r\^{o}le of the model
setting is one of clarification, the main question of this paper
can nevertheless be tackled in $\mathbb{F}_2^n$, and in particular
the following, for example, is proved in \cite{TS.L1NFTCVS}.
\begin{theorem}
Suppose that $G=\mathbb{F}_2^n$ and $A \subset G$ has density
$\alpha$ with $|\alpha - 1/3| \leq \epsilon$. Then
\begin{equation*}
\|\chi_A\|_{A(G)} \gg \log \log \epsilon^{-1}.
\end{equation*}
\end{theorem}

Finally it would be interesting to know what the true bounds in
Theorem \ref{mainconttheorem} should be. As far as the model
analogue, Theorem \ref{f2ndetail}, is concerned it would probably
be surprising if one could beat the following.
\begin{conjecture}
Suppose that $G=\mathbb{F}_2^n$, $f \in A(G)$ and $\epsilon \in
(0,1]$. Write $A_f$ for the quantity
$\|f\|_{A(G)}\|f\|_{L^\infty(G)}^{-1}$. Then there is a subspace $V$
of $G$ with
\begin{equation*}
\codim V \ll \epsilon^{-2} A_f,
\end{equation*}
and
\begin{equation*}
\sup_{x' \in G}{\|f-f(x')\|_{L^2(x'+V)}} \leq
\epsilon\|f\|_{L^\infty(G)}.
\end{equation*}
\end{conjecture}
It is, however, not clear what an argument giving this might provide
in the general setting. If the argument is iterative in the style of
this paper then to provide an improvement in the exponent of $\log
p$ in Theorem \ref{maintheorem} one would require some way of
cutting down the number of times we iterate.

\section*{Acknowledgments}
I should like to thank Tim Gowers for reviewing many drafts of the
paper, Ben Green for guidance and innumerable valuable
conversations, Ben Green and Sergei Konyagin for making some early
drafts of the preprint \cite{BJGSVK} available, and an anonymous
referee for a number of useful suggestions and improvements.

\appendix

\section{The general construction of
the auxiliary measure}\label{app1}

We extend the approach of \S\ref{auxiliarymeasure} from the model
setting to that of $G$ an arbitrary compact Abelian group. Here
Riesz products are marginally more complicated.

\subsection{Riesz products}\label{RPs2}

Suppose that $\Lambda$ is a finite set of characters. We say that
$\omega \in \ell^{\infty}(\Lambda \cup \Lambda^{-1})$ is
\emph{hermitian} if
\begin{equation*}
\omega(\lambda^{-1}) = \overline{\omega(\lambda)} \textrm{ for all
} \lambda \in \Lambda;
\end{equation*}
if $\omega$ also satisfies $\|\omega\|_\infty \leq 1$ then we
define the product
\begin{equation}\label{rieszproduct2}
p_\omega:=\prod_{\lambda \in
\Lambda}{\left(1+\frac{\omega(\lambda)\lambda +
\omega(\lambda^{-1})\lambda^{-1}}{2}\right)}.
\end{equation}
Such a product is called a \emph{Riesz product} and it is easy to
see that it is real and non-negative from which it follows that
$\|p_{\omega}\|_1 = \wh{p_{\omega}}(0_{\wh{G}})$. Further expanding
out the product reveals that $\supp \wh{p_\omega} \subset \langle
\Lambda \rangle$.

We had an easy time computing the Fourier transform of Riesz
products in $\mathbb{F}_2^n$, in general it is more complicated.
We can expand out the product in (\ref{rieszproduct2}) to see
that:
\begin{equation}\label{rpft}
\wh{p_{\omega}}(\gamma)=\sum_{m:\Lambda^m=\gamma}{\prod_{\lambda \in
\Lambda: \atop{m_\lambda \neq
0}}{\frac{\omega(\lambda^{m_{\lambda}})}{2}}}.
\end{equation}
To keep track of this we say that $\wt{p}$, defined on
$\{-1,0,1\}^{\Lambda}$, is a \emph{formal Fourier
transform}\footnote{Formal Fourier transforms are not in general
unique. } for $p \in L^1(G)$ if
\begin{equation}\label{fftfulfil}
\wh{p}(\gamma)=\sum_{m:\Lambda^m=\gamma}{\wt{p}(m)} \textrm{ for
all } \gamma \in \wh{G}.
\end{equation}
The functions which we are interested in are of the form
\begin{equation*}
p(x):=\int{p_{t\omega}(x)d\tau(t)},
\end{equation*}
for $\omega \in \ell^\infty(\Lambda \cup \Lambda^{-1})$ hermitian
with $\|\omega\|_\infty \leq 1$, and $\tau$ a real measure on
$[-1,1]$. It follows from (\ref{rpft}) and linearity of the
Fourier transform that $\wt{p}$ defined by
\begin{equation}\label{fft12}
\wt{p}(m):=\int{t^{|m|}d\tau(t)}.\prod_{\lambda \in \Lambda:
\atop{m_\lambda \neq 0}}{\frac{\omega(\lambda^{m_{\lambda}})}{2}}
\textrm{ for all } m \in \{-1,0,1\}^{\Lambda},
\end{equation}
is a formal Fourier transform for $p$.

If $\Lambda$ is dissociated then when $\gamma=0_{\wh{G}}$ there is
only one summand in the expression for $\wh{p_\omega}(\gamma)$
given in (\ref{rpft}) and that has a value of 1, so by
non-negativity of $p_\omega$
\begin{equation*}
\|p_\omega\|_1=\wh{p_\omega}(0_{\wh{G}})=1.
\end{equation*}
Dissociativity makes computing the Fourier transform easy for
$\gamma=0_{\wh{G}}$ by restricting the number of non-zero summands
in (\ref{rpft}); a lemma of Rider \cite{DR} provides a result for
more general $\gamma$:
\begin{lemma}\label{rider}
Suppose that $G$ is a compact Abelian group and $\Lambda$ is a
finite dissociated set of characters on $G$.\footnote{The
definition of dissociativity did not require any topological
structure on the group and indeed Rider's lemma is true without
the assumption that $\Lambda$ is a set of characters for simple
reasons of duality: if $H$ is the group generated by $\Lambda$
then we can endow it with the discrete topology and embed
$\Lambda$ in $\wh{\wh{H}}$. The image of $\Lambda$ under this
embedding is a dissociated set of characters on the compact
Abelian group $\wh{H}$, so that there is no real loss of
generality in assuming that $\Lambda$ is a set of characters.}
Then for all $\gamma \in \wh{G}$
\begin{equation*}
|\{m \in \{-1,0,1\}^{\Lambda}:|m|=r,\Lambda^m=\gamma\}| \leq 2^r.
\end{equation*}
\end{lemma}
\begin{proof}
Let $\omega$ be the hermitian function which takes $\Lambda$ to 1.
For this choice of $\omega$ (\ref{rpft}) is
\begin{equation*}
\wh{p_\omega}(\gamma)=\sum_{r \geq 0}{2^{-r}|\{m \in
\{-1,0,1\}^{\Lambda}:|m|=r,\Lambda^m=\gamma\}| }
\end{equation*}
But $|\wh{p_\omega}(\gamma)| \leq \|p_\omega\|_1=1$ since
$\Lambda$ is dissociated which yields the conclusion.
\end{proof}
\begin{proposition}\label{nearly}
Suppose that $G$ is a compact Abelian group, $\Lambda$ a finite
dissociated set of characters on $G$ with no elements of order 2 and
$\omega \in \ell^\infty(\Lambda \cup \Lambda^{-1})$ is hermitian
with $\|\omega\|_\infty \leq 1$. Then for any $\eta \in (0,1]$ there
is a function $f_\eta \in L^1(G)$ such that
\begin{equation*}
\wh{f_\eta}|_{\Lambda\cup\Lambda^{-1}}=\omega, \|f_{\eta}\|_1 \ll
(1+\log \eta^{-1}) \textrm{ and } |\wh{f_{\eta}}(\gamma)| \leq \eta
\textrm{ for all } \gamma \not \in \Lambda \cup \Lambda^{-1}.
\end{equation*}
\end{proposition}
\begin{proof}
Fix an integer $l>1$ to be optimized later and let $\tau_{2l}$ be
the measure of Lemma \ref{eqnlop}. Define
\begin{equation*}
p(x):=\int{p_{t\omega}d\tau_{2l}(t)},
\end{equation*}
and let $\wt{p}$ be the formal Fourier transform for $p$ defined by
(\ref{fft12}). $\wt{p}(m)=0$ if $|m|=0$ by definition of $\tau_{2l}$
and $\wt{p}$, so
\begin{eqnarray}
\nonumber \left|\wh{p}(\gamma) - \sum_{ |m|= 1
\atop{\Lambda^m=\gamma}}{\wt{p}(m)}\right| & \leq & \sum_{r \geq
2}{\sum_{ |m|= r \atop{\Lambda^m=\gamma}}{|\wt{p}(m)}|} \textrm{
by
definition (\ref{fftfulfil}),}\\
\nonumber & \leq & \sum_{r > 2l}{\sum_{ |m|= r
\atop{\Lambda^m=\gamma}}{|\wt{p}(m)}|} \textrm{ since $\int{t^rd\tau_{2l}(t)}=0$ for $1<r \leq 2l$,}\\
\nonumber & \leq & \sum_{r > 2l}{\sup_{|m|= r
\atop{\Lambda^m=\gamma}}{|\wt{p}(m)|}|\{m:
|m|=r, \Lambda^m=\gamma\}|}\\
\nonumber & \leq & \sum_{r > 2l}{2^r\sup_{ |m|= r
\atop{\Lambda^m=\gamma}}{|\wt{p}(m)|}} \textrm{ by Lemma
\ref{rider},}\\ \nonumber & \leq & \sum_{r >
2l}{2^r\int{t^rd\tau_{2l}(t)}(2^{-1}\|\omega\|_\infty)^r}\\
\nonumber & \leq & \sum_{r >
2l}{2(2^{-1}\|\omega\|_\infty)^r} \textrm{ since $|\int{t^rd\tau_{2l}(t)}|\leq 2^{1-r}$,}\\
\label{1bd} & \leq & 2^{1-2l}\|\omega\|_{\infty}.
\end{eqnarray}
Lemma \ref{rider} applies above because $\Lambda$ is dissociated.
Now let $l$ be such that $2^{3-2l} \leq \eta$ but $l \ll (1+\log
\eta^{-1})$ and put $f_\eta^{(1)}:=2p$. Then
\begin{enumerate}
\item If $\gamma \in \Lambda \cup \Lambda^{-1}$ then
\begin{equation*}
\sum_{|m|=1\atop{\Lambda^m=\gamma}}{\wt{p}(m)} =
\int{td\tau_{2l}(t)}\frac{\omega(\gamma)}{2} =
\frac{\omega(\gamma)}{2}
\end{equation*}
since $\Lambda$ has no elements of order 2. Hence by (\ref{1bd})
\begin{equation}
\label{1bound} |\wh{f_\eta^{(1)}}(\gamma)-\omega(\gamma)| \leq
2^{2-2l}\|\omega\|_{\infty} \leq 2^{-1}\|\omega\|_{\infty}.
\end{equation}
 \item If $\gamma \not\in \Lambda\cup \Lambda^{-1}$ then
\begin{equation*} \sum_{|m|=1 \atop{
\Lambda^m=\gamma}}{\wt{p}(m)} =0,
\end{equation*}
so by (\ref{1bd})
\begin{equation}
\label{otherbound} |\wh{f_{\eta}^{(1)}}(\gamma)| \leq
2^{-1}\eta\|\omega\|_{\infty}.
\end{equation}
\item $\|f_{\eta}^{(1)}\| \leq 2\|\tau_{2l}\|$ by the definition of $p$ and the
triangle inequality. \item $\wh{f_{\eta}^{(1)}}|_{\Lambda \cup
\Lambda^{-1}}$ is hermitian since $\tau_{2l}$ is real.
\end{enumerate}
We apply the foregoing recursively to the hermitian functions
$\omega$, $2(\omega - \wh{f_{\eta}^{(1)}})|_{\Lambda \cup
\Lambda^{-1}}$, $2(2(\omega - \wh{f_{\eta}^{(1)}}) -
\wh{f_{\eta}^{(2)}}) |_{\Lambda \cup \Lambda^{-1}}$,... to get a
sequence of $L^1(G)$-functions $f_{\eta}^{(1)}$, $f_{\eta}^{(2)}$,
$f_{\eta}^{(3)}$ ,... such that:
\begin{enumerate}
\item If $\gamma \in \Lambda \cup \Lambda^{-1}$ then
\begin{equation*}
|\sum_{k=1}^n{2^{-(k-1)}\wh{f_{\eta}^{(k)}}(\gamma)} -
\omega(\gamma)| \leq 2^{-n}.
\end{equation*}
\item If $\gamma \not \in \Lambda \cup \Lambda^{-1}$ then
\begin{equation*}
|\sum_{k=1}^n{2^{-(k-1)}\wh{f_\eta^{(k)}}(\gamma)}| \leq
\sum_{k=1}^n{2^{-(k-1)}.\frac{\eta}{2}} \leq \eta.
\end{equation*}
\item
\begin{equation*}
\|\sum_{k=1}^n{2^{-(k-1)}f_\eta^{(k)}}\|_1 \leq
\sum_{k=1}^n{2^{-(k-1)}\|f_\eta^{(k)}\|_1} \leq 2^2\|\tau_{2l}\|
\end{equation*}
\end{enumerate}
The sum $\sum_{k=1}^n{2^{-(k-1)}f_\eta^{(k)}}$ converges to a
function $f_\eta \in L^1(G)$ with the required properties since
$\|\tau_{2l}\| \ll l \ll (1+\log \eta^{-1})$.
\end{proof}
Finally we modify the above proposition so that the Fourier
transform is small on $\Lambda^{-1}\setminus \Lambda$.
\begin{proof}[Proof of Proposition \ref{auxiliary}] Let $G':=G \times S^1$ and
identify its dual with $\wh{G} \times \mathbb{Z}$; let
$\Lambda'=\Lambda \times \{1\}$, which is dissociated since
$\Lambda$ is dissociated, and has no elements of order 2 since 1
is not of order 2 in $\mathbb{Z}$; let $\omega'$ be the hermitian
map on $\Lambda' \cup \Lambda'^{-1}$ induced by
$\omega'(\lambda,1):=\omega(\lambda)$. Apply Proposition
\ref{nearly} to $G'$, $\Lambda'$ and $\omega'$ to get the function
$f_\eta \in L^1(G')$. Let $\mu_\eta$ be the measure induced by the
functional
\begin{equation*}
f \mapsto \int_{(x,z) \in
G'}{f(x)f_\eta(x,z)\overline{z}d\mu_{G'}(x,z)}
\end{equation*}
on $C(G)$. If $\gamma \in \wh{G}$ then
\begin{equation*}
\wh{\mu_\eta}(\gamma)=\int_{(x,z) \in
G'}{\overline{\gamma}(x)f_\eta(x,z)\overline{z}d\mu_{G'}(x,z)}=\wh{f_\eta}(\gamma,1).
\end{equation*}
We verify the three properties of $\mu_\eta$ from the
corresponding properties of $f_\eta$:
\begin{enumerate}
\item If $\lambda \in \Lambda$ then
$\wh{\mu_\eta}(\lambda)=\wh{f_\eta}(\lambda,1)=\omega'(\lambda,1)=\omega(\lambda)$.
\item
\begin{eqnarray*}
\|\mu_\eta\| & = & \sup_{f \in C(G):\|f\|_\infty \leq
1}{\left|\int_{(x,z) \in
G'}{\overline{\gamma}(x)f_\eta(x,z)\overline{z}d\mu_{G'}(x,z)}\right|}\\
& \leq & \int_{(x,z) \in G'}{|f_\eta(x,z)|d\mu_{G'}}\\ & \ll &
(1+\log \eta^{-1}).
\end{eqnarray*}
\item If $\gamma \not \in \Lambda$ then $(\gamma,1) \not \in
\Lambda' \cup \Lambda'^{-1}$ so $|\wh{\mu_\eta}(\gamma)| \leq
\eta$.
\end{enumerate}
\end{proof}

\bibliographystyle{alpha}

\bibliography{master}

\end{document}